\documentclass[11pt,reqno]{amsart}

\usepackage{hyperref}

\allowdisplaybreaks[4]

\usepackage{amsmath}
\usepackage{amssymb}
\usepackage{mathrsfs}
\usepackage{amsthm}
\usepackage{bm}
\usepackage{graphicx}
\usepackage{booktabs}
\usepackage{float}
\usepackage{geometry}
\usepackage{color}

\theoremstyle{plain}
\newtheorem{lemma}{Lemma}[section]
\newtheorem{theorem}[lemma]{Theorem}

\newtheorem{remark}[lemma]{Remark}
\newtheorem{proposition}[lemma]{Proposition}

\newtheorem{assumption}{Assumption}

\newcommand{\N}{\mathbb{N}}
\newcommand{\R}{\mathbb{R}}
\newcommand{\E}{\mathbb{E}}

\newcommand{\IW}{\mathcal{W}}

\newcommand{\IF}{\mathcal{F}}

\newcommand{\IL}{\mathcal L}

\newcommand{\<}{\langle}

\begin{document}
\title[$L^p$-strong convergence orders of fully discrete schemes]{$L^p$-strong convergence orders of fully discrete schemes for  the SPDE driven by L\'evy noise}
\author{Chuchu Chen, Tonghe Dang, Jialin Hong, Ziyi Lei}
\address{LSEC, ICMSEC,  Academy of Mathematics and Systems Science, Chinese Academy of Sciences, Beijing 100190, China,
\and 
School of Mathematical Sciences, University of Chinese Academy of Sciences, Beijing 100049, China,
\and Department of Applied Mathematics, The Hong Kong Polytechnic University, Hung Hom, Kowloon, Hong Kong}
\email{chenchuchu@lsec.cc.ac.cn; tonghe.dang@polyu.edu.hk; hjl@lsec.cc.ac.cn; ziyilei@lsec.cc.ac.cn}

\thanks{This work is funded by the National key R\&D Program of China under Grant (No. 2024YFA1015900 and No. 2020YFA0713701), National Natural Science Foundation of China (No. 12031020, No. 12461160278, No. 12471386,  and No. 12288201), and by Youth Innovation Promotion Association CAS, China.}
\begin{abstract}  
It is well known that for a stochastic differential equation driven by L\'evy noise, the temporal H\"older continuity in $L^p$ sense of the exact solution does not exceed $1/p$. This leads to that the $L^p$-strong convergence order of a numerical scheme will vanish as $p$ increases to infinity if the temporal H\"older continuity of the solution process is directly used. A natural question arises: can one obtain the $L^p$-strong convergence order that does not depend on $p$? In this paper, we provide a positive answer for fully discrete schemes of the stochastic partial differential equation (SPDE) driven by L\'evy noise. Two cases are considered: the first is the linear multiplicative Poisson noise with $\nu(\chi)<\infty$ and the second is the additive Poisson noise with $\nu(\chi)\leq\infty$, where $\nu$ is the L\'evy measure and $\chi$ is the mark set. For the first case, we present a strategy by employing the jump-adapted time discretization, while for the second case, we introduce the approach based on the recently obtained L\^e's quantitative John--Nirenberg inequality. We show that proposed schemes converge in $L^p$ sense with orders almost $1/2$ in both space and time for all $p\ge2$, which contributes novel results in the numerical analysis of the SPDE driven by L\'evy noise. 
  \end{abstract}
\keywords {Fully discrete scheme $\cdot$ $L^p$-strong convergence order $\cdot$ L\'evy noise $\cdot$ Stochastic partial differential equation}
\maketitle

\section{Introduction}

Stochastic differential equations (SDEs) with L\'evy noise are widely used to model sudden events and irregular changes in stochastic phenomenon, such as financial market crashes, abrupt phase transitions, and neural spiking patterns; see e.g. \cite{RP2004,MR2356959,W2001} and references therein. In numerical studies of SDEs, the $L^p$-strong convergence analysis of numerical methods has received much attention, as it assesses the accuracy of simulating the paths of the underlying solution process. The $L^p$-strong convergence analysis for SDEs with L\'evy noise is different from that of the Gaussian noise case (see e.g. \cite{MR3022246,MR4408264}). For instance, there is an extra term appearing in the Burkholder--Davis--Gundy inequality of the stochastic integral with respect to the compensated Poisson measure compared with that in the Gaussian noise case. As a result, the temporal H\"older continuity in $L^p$ sense of the exact solution to an SDE with L\'evy noise does not exceed $1/p$, resulting in the deterioration of the $L^p$-strong convergence order of a numerical scheme for large $p$ if one uses the temporal H\"older continuity directly. A natural question arises: 
\begin{itemize}
\item[Q:] \textit{Can one obtain the $L^p$-strong convergence order of a numerical scheme for the SDE with L\'evy noise that does not depend on $p$? }
\end{itemize}

A remarkable progress for this question has been made recently for the case of the stochastic ordinary differential equation (SODE) with L\'evy noise: authors in \cite{OKM-2024} studied the $L^p$-strong convergence order of the Euler--Maruyama scheme for a multidimensional SODE with irregular H\"older drift, driven by additive L\'evy process with exponent $\alpha\in(0,2]$, and the obtained convergence order does not depend on $p$. For the case of SPDEs, to the best of our knowledge, this question is still open so far. The aim of this paper is to provide a positive answer for fully discrete schemes of the SPDE driven by L\'evy noise.

We consider the following SPDE driven by L\'evy noise:
\begin{align}\label{SPDE:MultiplicativeNoise}
\begin{cases}
dX(t)=AX(t)dt+F(X(t))dt+dW(t)+\int_{\chi}G(X(t), z)\tilde{N}(dz,dt),\ \ t\in(0,T],\\
X(0)=x_0,
\end{cases}
\end{align}
where $T>0$ and $A:=\Delta:\mbox{Dom}(A)\subset H\to H$ is the Laplacian with homogeneous Dirichlet boundary condition. Here $H:=L^2(0,1)$ with usual inner product $\langle\cdot,\cdot\rangle$ and norm $\|\cdot\|$, and $\chi:=H\backslash\{0\}$ is the mark set. The process $\{W(t)\}_{t\in[0,T]}$ is an $H$-valued Gaussian space-time white noise defined on a complete filtered probability space $(\Omega, \IF,(\IF_t)_{t\geq0},\mathbb P)$ that satisfies usual conditions. Let $\tilde{N}(dz,dt):=N(dz,dt)-\nu(dz)dt$ stand for a compensated Poisson random measure, where $\{N(\cdot,t)\}_{t\in[0,T]}$ is the Poisson random measure and $\nu$ is a L\'evy measure on Borel $\sigma$-algebra $\mathcal B(\chi)$ satisfying $\nu(\{0\})=0$ and $\int_{\chi}\min\{1,\|z\|^2\}\nu(dz)<\infty.$ Here $\{W(t)\}_{t\in[0,T]}$ and $\{\tilde{N}(\cdot,t)\}_{t\in[0,T]}$ are supposed to be independent. Precise assumptions on coefficients $F$ and $G$ will be given in Section \ref{sec2}. For the convergence analysis of numerical methods for ~\eqref{SPDE:MultiplicativeNoise}, we refer interested readers to \cite{CTJ2024} for the $L^2$-strong convergence and to \cite{E2007} for the $L^p$-strong convergence with $p\leq 2$.

 We apply the spectral Galerkin method to approximate  \eqref{SPDE:MultiplicativeNoise} in space, and further use the Euler-type methods in time to obtain the fully discrete schemes. The corresponding mapping for the one-step approximation is given as follows: 
\begin{align*}
		\Phi(X^N_{i},\Delta t_{ i+1})&:=E(\Delta t_{ i+1})X^N_{i}+\int_{t_i}^{t_{i+1}}E(t_{i+1}-s)F_N(X^N_{i})ds\nonumber\\
		&~\quad+\int_{t_i}^{t_{i+1}}E(t_{i+1}-s)P_NdW(s)\nonumber\\
		&~\quad+\int_{t_i}^{t_{i+1}}\!\int_\chi E(\Delta t_{ i+1})G_N(X^N_{i}, z)\tilde{N}(dz,ds),\quad i=0,1,\ldots,n_T-1,
\end{align*}
where $n_T$ is the subscript corresponding to the time node $T$ (i.e., $t_{n_T}=T$), and we mention that the time step size $\Delta t_{ i+1}:=t_{i+1}-t_i$ may be path-dependent. The main contribution of our work is that, we show for the first time that proposed schemes converge in $L^p$ sense with orders almost $1/2$ in both space and time for all $p\ge2$. Indeed, the $L^p$-strong convergence analysis in time is more technical compared with that in space. This is related to the aforementioned  issue yielding poor order for large $p$. To overcome this order barrier, we present two different strategies for the cases of the linear multiplicative Poisson noise with $\nu(\chi)<\infty$ and the additive Poisson noise with $\nu(\chi)\leq\infty$, respectively. More precisely,

 \begin{itemize}
\item[(\romannumeral1)]
 \textit{Case of the linear multiplicative Poisson noise with $\nu(\chi)<\infty$.} 
In this case, there are a.s. finitely many jumps.  We introduce the jump-adapted time discretization by a superposition of finitely many jump times to a deterministic equidistant grid. Then the jump-adapted fully discrete scheme is constructed  as follows: 
   \begin{align}\label{eq:fullyscheme}
     X^N_{i+1-}=\Phi(X^N_i,\Delta t_{i+1}), \quad X^N_{i+1}=X^N_{i+1-}+G_N(X^N_{i+1-},p_{t_{i+1}})\mathbf{1}_{\{p_{t_{i+1}}\neq 0\}}. 
   \end{align}
Here $\{p_{t_i}\}_{i=0,1,..., n_T}$ is an $\mathcal F_{t_i}$-adapted Poisson point process with the intensity measure $\nu(dz)dt$. In the convergence analysis of this scheme, there are two critical steps to overcome the order barrier. First, this scheme evolves without jump between time grid points, which implies that the stochastic integral with respect to the compensated Poisson random measure  in~\eqref{eq:fullyscheme} can be transformed into the integral with respect to the intensity measure. Second, it requires establishing the $p$th moment estimates on the terms involving the multiplicative noise or the jump process when accumulating the local errors into the global one.

\item[(\romannumeral2)]
 \textit{Case of the additive Poisson noise with $\nu(\chi)\leq\infty$.} In this case, there may be a.s. infinitely many jumps, in which the jump-adapted strategy is invalid. We consider the fully discrete scheme $X^N_{i+1}=\Phi(X^N_i,\Delta t)$ with constant step size $\Delta t$. In the $L^p$-strong convergence analysis of this scheme, new strategy is required to deal with the terms involving the temporal H\"older continuity of the stochastic convolution with respect to the compensated Poisson random measure of \eqref{SPDE:MultiplicativeNoise}. Our idea is to bound the $L^p$-norm of these terms by moment estimates of the nested conditional ``$L^1$-norm"  based on the recently obtained L\^e's quantitative John--Nirenberg inequality.
To this end, we carefully analyze some essential \textit{a priori} estimates including the nested conditional ``$L^2$-norms" of both the stochastic convolutions and the solution of the perturbed stochastic equation, which help us to obtain the rigorous bound of the nested conditional ``$L^1$-norm" and then the desired convergence orders. 
\end{itemize}

The rest of the paper is organized as follows. The preliminaries, including notations, assumptions, and the well-posedness of the exact solution are given in Section \ref{sec2}. Section \ref{subsec2.2} is devoted to the introduction of the fully discrete schemes and their $L^p$-strong convergence order results. In Section \ref{sec3}, we present the $L^p$-strong convergence analysis of the fully discrete scheme for \eqref{SPDE:MultiplicativeNoise} in the case of the linear multiplicative Poisson noise with $\nu(\chi)<\infty$.  In Section \ref{sec4}, we present the corresponding $L^p$-strong convergence analysis for the case of the additive Poisson noise with $\nu(\chi)\leq\infty$. Proofs of some essential propositions used in the convergence analysis are given in Section \ref{sec5}.

\section{Fully discrete schemes and main results}

In this section, we first give some preliminaries, including notations, assumptions, the well-posedness of the exact solution and some useful inequalities. Second, we present the fully discrete schemes of \eqref{SPDE:MultiplicativeNoise} for two different cases of noises, and respectively show their $L^p$-strong convergence orders, as the main results of this paper.

\subsection{Preliminaries}\label{sec2}

Throughout this paper, we let $C>0$ be a generic constant that may vary from one place to another. More specific constants which depend on certain parameters $a,b$ are numbered as $C_{a,b}$. We use $\N_+$ to denote the set of positive natural numbers and let $\epsilon,\delta>0$ be arbitrarily small parameters. Denote $a\vee b:=\max\{a,b\}$. The space of bounded linear operators on $H$ is denoted by $\IL(H)$  with the operator norm $\|\cdot\|_{\IL(H)}$. The subspace $\IL_2(H)\subset\IL(H)$ denotes the set of Hilbert--Schmidt operators on $H$, with norm denoted by $\|\cdot\|_{\IL_2(H)}$. Let $\dot{H}^s$, $s\in\R$ be the Sobolev space generated by the fractional power of $-A$, endowed with the inner product $\langle u,v\rangle_{\dot{H}^s}:=\langle(-A)^{\frac{s}{2}}u,(-A)^{\frac{s}{2}}v\rangle$ and the norm $\|u\|_{\dot{H}^s}:=\sqrt{\langle u,u\rangle_{\dot{H}^s}}$ for all $u,v\in \dot{H}^s$. We use the notation $L^p(\Omega,\dot{H}^{s})$ to denote the space of random variables $u$ satisfying $\E[\|u\|^p_{\dot{H}^s}]<\infty$ for $p>0$.

It is known that there is a sequence of real numbers $\lambda_i=\pi^2i^2, i\in\N_+$ and an orthonormal basis $\{e_i\}_{i\in\N_+}$ of $H$ with $e_i(x):=\sqrt{2}\sin(i\pi x)$, $x\in[0,1]$, such that $-Ae_i=\lambda_ie_i$. In addition, operator $A$ generates the semigroup $\{E(t):=e^{tA}, t\geq0\}$ on $H$. For all $t\geq0$, $E(t)$ is a bounded self-adjoint linear operator on $H$, with $\|E(t)\|_{\IL(H)}\le e^{-t\lambda_1}\leq 1$. It is clear that, for all $t\ge 0$, 
   \[\int_0^ts^{-\alpha}\|E(s)\|_{\IL_2(H)}^2ds<\infty\quad\quad\mbox{if and only if}\quad\alpha\in[0,\frac12).
   \]
Moreover, the semigroup satisfies the following properties: 
\begin{align}
\|(-A)^{\gamma}E(t)\|_{\IL(H)}&\leq Ct^{-\gamma}e^{-\frac{\lambda_1}{2}t},\quad t\in(0,T],~\gamma\ge 0,\label{prop:semigroup1}\\
	\|(-A)^{-\rho}(E(t)-E(s))\|_{\IL(H)}&\leq C(t-s)^{\rho},\quad 0\leq s<t\leq T, ~\rho\in[0,1].\label{prop:semigroup2}
\end{align}

Here we impose some assumptions on coefficients and the initial value of \eqref{SPDE:MultiplicativeNoise}, which will be used throughout this paper.

\begin{assumption}[Nonlinearity]\label{assum:nonlinearity}
  The measurable mapping $F\colon H\to H$ satisfies that $\|F(x)\|\leq C(1+\|x\|)$, $x\in H$ with some constant $C>0$. In addition, $F$ is differentiable and there are constants $L_F,C>0$ such that
 \begin{align*}
	\sup_{x\in H}\|DF(x)y\|&\leq L_F\|y\|,\quad y\in H,\\
     \|(-A)^{-\frac{1+\delta}{4}}DF(x)y\|&\leq C(1+\|x\|_{\dot{H}^{\frac{1-\delta}{2}}})\|y\|_{\dot{H}^{-\frac{1-\delta}{2}}},\quad x\in\dot{H}^{\frac{1-\delta}{2}}, y\in\dot{H}^{-\frac{1-\delta}{2}}.
\end{align*} 
\end{assumption}

\begin{remark}
   The conditions in Assumption~\ref{assum:nonlinearity} are satisfied if $F$ is the Nemytskii operator defined by
$F(u)(x):=f(u(x)),\; x\in[0,1],\;u\in H,$ where $f\colon \R\to\R$ is continuously differentiable with bounded first order derivative. Indeed, the mapping $F\colon H\to H$ is then differentiable, and for all $u,h\in H$ and $x\in[0,1]$, we have
$
       [DF(u)h](x)=f'(u(x))h(x).
$
This implies $\sup_{u\in H}\|DF(u)h\|\leq L_F\|h\|, h\in H$ with some constant $L_F>0.$ In addition, 
for all $u_1\in \dot{H}^{\frac{1-\delta}{2}}$ and $u_2\in\dot{H}^{-\frac{1-\delta}{2}}$, we obtain
\begin{align*}
   \big\|(-A)^{-\frac{1+\delta}{4}}DF(u_1)u_2\big\|&=\sup_{h\in H\setminus\{0\}}\frac{\<(-A)^{-\frac{1+\delta}{4}}DF(u_1)u_2,h\rangle}{\|h\|}\\
    &\leq \sup_{h\in H\setminus\{0\}}\frac{\|(-A)^{-\frac{1-\delta}{4}}u_2\|\|(-A)^{\frac{1-\delta}{4}}DF(u_1)(-A)^{-\frac{1+\delta}{4}}h\|}{\|h\|}\\
    &\leq C(1+\|(-A)^{\frac{1-\delta}{4}}u_1\|)\|(-A)^{-\frac{1-\delta}{4}}u_2\|.
\end{align*}
where we used \cite[Lemma 4.4]{XGT2014} in the last step. 
\end{remark}

\begin{assumption}[Jump coefficient]\label{assum:jumpcoefficient}
	The coefficient $G\colon H\times\chi\to H$ is defined as $G(x,z)=g_1(z)x+g(z),\,x\in H,\,z\in\chi$ with mapping $g_1\colon \chi\to\R$ being bounded by constant $b>0$ and with some mapping $g\colon \chi\to H$. In addition, for each $p>2,$ there exists a constant $C_{p}>0$ such that
\begin{align*}
	\int_\chi\big(|g_1(z)|^2\vee|g_1(z)|^p\big)\nu(dz)\leq C_p,\quad
	\int_\chi \big(\|(-A)^{\frac14}g(z)\|^2\vee\|(-A)^{\frac14}g(z)\|^p\big)\nu(dz)\leq C_{p}.
	\end{align*} 
\end{assumption}

\begin{assumption}[Initial value]\label{assum:x0}
Let $x_0$ be an $\IF_{0}$-measurable random variable,  and $x_0\in L^{p}(\Omega,\dot{H}^{\frac{1-\delta}{2}})$ for all $p\ge 2$.
\end{assumption}

We introduce the maximal inequality for the L\'evy-type stochastic convolution as follows, see e.g. \cite[Proposition 3.3]{CCM2010} for details.
\begin{lemma}\label{ineq:maximalinequality}
 Let $\beta:[0,T]\times \chi\to H$ be a predictable process such that the expectation on the right hand side of the inequality below is finite. Then for all $p\ge 2$, there exists a constant $C_{p,T}>0$ such that
\begin{align*}
    &\quad\E\Big[\sup_{t\in[0,T]}\Big\|\int_0^t\int_\chi E(t-s)\beta(s,z)\tilde{N}(dz,ds)\Big\|^p\Big]\\
	&\leq C_{p,T}\E\Big[\int_0^T\Big\{\Big(\int_\chi\|\beta(s,z)\|^2\nu(dz)\Big)^{p/2}+\int_\chi\|\beta(s,z)\|^p\nu(dz)\Big\}ds \Big].
\end{align*}   
\end{lemma}

Under assumptions given above, the mild solution of \eqref{SPDE:MultiplicativeNoise} uniquely exists and the regularity estimate of the solution can be obtained, which are stated in the following proposition. The proof of the well-posedness is similar to that of \cite[Theorem 2.1]{BX-2012} and the one of the strong Markov property is similar to that of \cite[Lemma 5.2]{BX-2012}, which are omitted here. For the proof of regularity estimate, we postpone it to Appendix.

\begin{proposition}\label{prop:wellposedness}
   Let Assumptions~\ref{assum:nonlinearity}--\ref{assum:x0} hold. Then for each $T>0,$  there exists a unique mild solution $\{X(t)\}_{t\in[0,T]}$ of  \eqref{SPDE:MultiplicativeNoise} which reads as
 \begin{align*}
X(t)&=E(t)x_0+\int_{0}^tE(t-s)F(X(s))ds+\int_{0}^tE(t-s)dW(s)\nonumber\\
   &\quad+\int_{0}^t\int_{\chi}E(t-s)G(X(s),z)\tilde{N}(dz,ds),\quad t\in(0,T]
  \end{align*} 
with the initial value $x_0$ and satisfies that for any $p\geq 2,$  
$
       \E[\sup_{t\in[0,T]}\|X(t)\|^p]<\infty.
$
  The process $\{X(t)\}_{t\in[0,T]}$ has a c\`{a}dl\`{a}g modification and is  time-homogeneous strong Markovian.  Moreover, for all $p\geq2$ and all $\alpha\in[0,\frac12)$, there exists a constant $C_{p,T}>0$ such that
\begin{align*}
\sup_{t\in[0,T]}\|X(t)\|_{L^p(\Omega, \dot{H}^\alpha)}\leq C_{p,T}.
	\end{align*} 
\end{proposition}

\subsection{Fully discrete schemes and $L^p$-strong convergence orders}\label{subsec2.2}

In this subsection, we first apply the spectral Galerkin method in the spatial direction and the Euler-type method in the temporal direction to obtain the fully discrete schemes for \eqref{SPDE:MultiplicativeNoise}. Then we present the main results of this paper, namely, $L^p$-strong convergence orders of proposed schemes in both space and time.

For fixed $N\in\N_+$, let $P_N$ be the orthogonal projection operator from $H$ onto $H_N:=\mbox{span}\{e_1,e_2,\ldots,e_N\}$. Applying the spectral Galerkin method to the spatial direction of \eqref{SPDE:MultiplicativeNoise} yields
\begin{align}\label{eq:spectralgalerkin}
dX^N(t)&=\big(AX^N(t)+F_N(X^N(t))\big)dt+P_N dW(t)\nonumber\\
		&\quad+\int_{\chi}G_N(X^N(t),z)\tilde{N}(dz,dt),~~ t\in(0,T]
\end{align}
with the initial value $X^N(0)=P_Nx_0,$ 
where $F_N:=P_NF$ and $G_N:=P_NG=g_{1}P_N+P_Ng$. Further, in the temporal direction, we use the Euler-type method to obtain the one-step approximation: 
 \begin{align}\label{EM_full}
		\Phi(X^N_{i},\Delta t_{ i+1})&=E(\Delta t_{ i+1})X^N_{i}+\int_{t_i}^{t_{i+1}}E(t_{i+1}-s)F_N(X^N_{i})ds\nonumber\\
		&~\quad+\int_{t_i}^{t_{i+1}}E(t_{i+1}-s)P_NdW(s)\nonumber\\
		&~\quad+\int_{t_i}^{t_{i+1}}\!\!\int_\chi E(\Delta t_{ i+1})G_N(X^N_{i}, z)\tilde{N}(dz,ds),\quad i=0,1,\ldots,n_T-1,
\end{align}
where $\Delta t_{i+1}=t_{i+1}-t_i$ with $\{t_i\}_{i=0,1,\ldots,n_T}$ being time grid points. With the one-step approximation at hand, below, we present the fully discrete schemes for two cases: the first is the linear multiplicative Poisson noise with $\nu(\chi)<\infty$, and the second is the additive Poisson noise with $\nu(\chi)\leq\infty$. We remark that the choices of time step sizes $\{\Delta t_{i}\}_{i=1,\ldots,n_T}$ are different in these two cases.

\textbf{(\romannumeral1) Case of the linear multiplicative Poisson noise with $\nu(\chi)<\infty$.}
In this case, since $\mathbb E[N(\chi,[0,T])]=T\nu(\chi)<\infty,$ there are a.s. finitely many jumps. The corresponding non-decreasing jump times are denoted by $\{\sigma_i(\omega),i=1,2,\ldots,n_J(\omega)\}$,  where the random variable $n_J$ denotes the number of jumps. We introduce the jump-adapted time partition by a superposition of finitely many jump times to a deterministic equidistant grid. More precisely, we first give a deterministic partition of the interval $[0,T]$:
$$
\mathcal{T}^0=\{0=t_{0}^0<t_{1}^0<\cdots<t_M^0=T\},
$$
where $\Delta t>0$ is a constant step size and $t_{i}^0:=i \Delta t, i=0, 1, \ldots, M$. For each sample path, we then merge $\mathcal{T}^0$ and the partition from jump times $\mathcal{T}^1:=\{0\leq \sigma_1<\sigma_2< \cdots<\sigma_{n_J}\leq T\}$  to form a new partition
$$
\mathcal{T}=\{0=t_{0}<t_{1}<\cdots<t_{n_{T}}=T\},
$$
where $n_{T}$ denotes the subscript corresponding to the time node $T$ and $n_J+M\ge n_T$. Note that in the new partition $\mathcal{T}$, step size $\Delta t_{i+1}=t_{i+1}-t_i$ is path-dependent, and the maximal time step size of the resulting jump-adapted time partition is no larger than $\Delta t$.

The jump-adapted fully discrete scheme is proposed as 
 \begin{align*}
     X^N_{i+1-}=\Phi(X^N_i,\Delta t_{i+1}), \quad X^N_{i+1}=X^N_{i+1-}+G_N(X^N_{i+1-},p_{t_{i+1}})\mathbf{1}_{\{p_{t_{i+1}}\neq 0\}}, 
   \end{align*}
with $\Phi(X^N_i,\Delta t_{i+1})$ being given by \eqref{EM_full}, which is equivalent to 
\begin{align}
\label{fullnumericalscheme}
\begin{cases}
X^N_{ i+1-}\!\!=E(\Delta t_{ i+1})X^N_{i}+\int_{t_i}^{t_{i+1}}\!E(t_{i+1}\!-\!s)F_N(X^N_{i})ds+\int_{t_i}^{t_{i+1}}\!E(t_{i+1}\!-\!s)P_NdW(s)\\
\quad\quad\qquad-\int_{t_i}^{t_{i+1}}\int_\chi E(\Delta t_{i+1})G_N(X^N_{i}, z)\nu(dz)ds,\\		X^N_{i+1}=X^N_{i+1-}+G_N(X^N_{i+1-},p_{t_{i+1}})\mathbf{1}_{\{p_{t_{i+1}}\neq 0\}}, \quad i=0,1,\ldots, n_T-1.
\end{cases}
	\end{align}
due to the following relation:
 \begin{align}\label{relation1}
	\int_{t_i}^{t_{i+1}}\!\!\!\int_\chi E(\Delta t_{ i+1})&G_N(X^N_{i}, z)\tilde{N}(dz,dt)\nonumber\\
	&=-\int_{t_i}^{t_{i+1}}\!\!\!\int_\chi  E(\Delta t_{ i+1})G_N(X^N_{i}, z) \nu(dz)dt,\quad a.s.
\end{align}
for $i=0,1,\ldots, n_T-1$. We now give the $L^p$-strong convergence result of the scheme \eqref{fullnumericalscheme}, whose proof is presented in Section \ref{sec3}.

\begin{theorem}\label{the:convergenceorder-jumpadapted}
Let Assumptions~\ref{assum:nonlinearity}--\ref{assum:x0} hold. Then for all $N\in\N_+$ and all $p\geq2$, there exists a constant $C_{p,T}>0$ independent of $N,\Delta t$ such that
  \begin{align*}
      \|X(T)-X^N_{n_{T}}\|_{L^p(\Omega,H)}\leq C_{p,T}\big(N^{-\frac{1-\delta}{2}}+(\Delta t)^{\frac{1-\delta}{2}}\big).
  \end{align*} 
\end{theorem}

\textbf{(\romannumeral2) Case of the additive Poisson noise with $\nu(\chi)\leq\infty$.} In this case when $\nu(\chi)<\infty$, we can obtain the same $L^p$-strong convergence orders as shown in Theorem \ref{the:convergenceorder-jumpadapted} by applying the scheme \eqref{fullnumericalscheme}. 
  However, this jump-adapted numerical scheme is invalid when $\nu(\chi)=\infty$, since there may be a.s. infinitely many jumps. In order to obtain the $L^p$-strong convergence orders in a wider case: $\nu(\chi)\leq\infty$, we consider the fully discrete scheme on the deterministic partition $\mathcal{T}^0$ with constant step size $\Delta t$. For simplicity, we write $t_{i}$ instead of time grid point $t_{i}^0, i=0, 1, \ldots, M$ and $n_T=M.$ Then the fully discrete scheme based on the one-step approximation~\eqref{EM_full} is given as 
\begin{align}\label{eq:fulldiscrete-onestep}
X^N_{k+1}=\Phi(X^N_k,\Delta t),\;k=0,1,\ldots,n_T-1,\quad X^N_0=P_Nx.
\end{align}
New strategy is proposed in the $L^p$-strong convergence analysis of this scheme. We now give the $L^p$-strong convergence result of the scheme \eqref{eq:fulldiscrete-onestep} for the case of the additive Poisson noise ($g_1\equiv0$) as follows, and its proof is presented in Section~\ref{sec4}.

\begin{theorem}\label{theo:temporalstrong-Lp}
Let Asumptions~\ref{assum:nonlinearity}--\ref{assum:x0} hold and $g_1\equiv0$. In addition, let $g$ satisfy that for all $p\geq2$,
\begin{align}\label{condition:g}
\int_{\chi}\big\|(-A)^{\frac{1-\delta}{2}}g(z)\big\|^p\nu(dz)<\infty.  
\end{align}
Then for all $N\in\N_+$ and all $p\geq2$, there exists a constant $C_{p,T}>0$ independent of $N,\Delta t$ such that
\begin{align*}
\big\|X(T)-X^N_{n_T}\big\|_{L^p(\Omega,H)}\leq C_{p,T}\big(N^{-\frac{1-\delta}{2}}+(\Delta t)^{\frac{1-\delta}{2}}\big).
\end{align*}
\end{theorem}

\section{ Proof of Theorem~\ref{the:convergenceorder-jumpadapted}}\label{sec3}
In this section, we present the $L^p$-strong convergence analysis of fully discrete scheme~\eqref{fullnumericalscheme} for \eqref{SPDE:MultiplicativeNoise} driven by the linear multiplicative Poisson noise with finite L\'evy measure. Note that  the stochastic integral with respect to the compensated Poisson random measure in \eqref{fullnumericalscheme} can be transformed into the integral with respect to the intensity measure, as shown in \eqref{relation1}. This transformation is crucial to overcome the order barrier caused by the $p$th moment estimates of the stochastic integral with respect to the compensated Poisson random measure.

\begin{proof}[\textbf{Proof of Theorem \ref{the:convergenceorder-jumpadapted}}] The proof is split into two steps, based on the error between  the solutions of \eqref{SPDE:MultiplicativeNoise} and the semi-discrete scheme \eqref{eq:spectralgalerkin}, and the error between the solutions of \eqref{eq:spectralgalerkin} and the fully discrete scheme \eqref{fullnumericalscheme}. 

\textbf{Step 1.} \textit{Show that for all $p\geq2$, there exists a constant $C_{p,T}>0$ independent of $N$ such that}
$
      \|X(T)-X^N(T)\|_{L^p(\Omega,H)}\leq C_{p,T}N^{-\frac{1-\delta}{2}}.
$  

 It follows from \eqref{SPDE:MultiplicativeNoise} and \eqref{eq:spectralgalerkin} that, for all $N\in\N_+$, 
  \begin{align*}
   X(T)-X^N(T)
    &=E(T)(\mathrm{I}-P_N)x_0+\int_0^T E(T\!-\!s)(\mathrm{I}-P_N)F(X(s))ds\\
    &\quad+\!\int_0^T E(T\!-\!s)P_N\big(F(X(s))-F(X^N(s))\big)ds\!+\!(\mathrm{I}-P_N)\int_0^TE(T\!-\!s)dW(s)\\
    &\quad+\!\int_0^T\int_\chi E(T\!-\!s)(\mathrm{I}-P_N)G(X(s),z)\tilde{N}(dz,ds)\\
    &\quad+\!\int_0^T\int_\chi E(T\!-\!s)P_N\big(G(X(s),z)-G(X^N(s),z)\big)\tilde{N}(dz,ds),
  \end{align*}
where $\mathrm{I}\colon H\to H$ is an identical mapping. Denote the stochastic convolution $\mathcal{W}(t):=\int_0^{t}E(t-s)dW(s)$. Owing to properties of projection operator $P_N$, assumptions on coefficients $F$ and $G$, Lemma~\ref{ineq:maximalinequality}, and Proposition~\ref{prop:wellposedness}, we obtain \begin{align*}
&\quad\E\big[\|X(T)-X^N(T)\|^p\big]\\
&\leq C_p\lambda_N^{-\frac{1-\delta}{4}p}\|x_0\|^p_{L^p(\Omega,\dot{H}^{\frac{1-\delta}{2}})}+C_p\lambda_N^{-\frac{1-\delta}{4}p}\Big(\int_0^T(T-s)^{-\frac{1-\delta}{4}}ds\Big)^p\sup_{s\in[0,T]}\E[\|F(X(s))\|^p]\\
&\quad+C_{p,T}\int_0^T\E\big[\big\|X(s)-X^N(s)\big\|^p\big]ds+C_p\lambda_N^{-\frac{1-\delta}{4}p}\big\|\IW(T)\|^p_{L^p(\Omega,\dot{H}^{\frac{1-\delta}{2}})}\\
 &\quad+C_{p,T}\lambda_N^{-\frac{1-\delta}{4}p}\E\Big[\int_0^T\!\!\!\Big\{\int_\chi\big\|(-A)^{\frac{1-\delta}{4}}G(X(s),z)\big\|^p\nu(dz)\!+\!\Big(\!\int_\chi\big\|(-A)^{\frac{1-\delta}{4}}G(X(s),z)\big\|^2\nu(dz)\Big)^{\frac{p}{2}}\Big\}ds\Big]\\
 &\leq C_{p,T}\lambda_N^{-\frac{1-\delta}{4}p}+C_{p,T}\int_0^T \!\!\E\big[\big\|X(s)-X^N(s)\big\|^p\big]ds\\
 &\quad+C_{p,T}\lambda_N^{-\frac{1-\delta}{4}p}\sup_{s\in[0,T]}\big\|X(s)\big\|^p_{L^p(\Omega,\dot{H}^{\frac{1-\delta}{2}})}.
\end{align*}
Applying the Gr\"onwall inequality finishes the proof of \textbf{Step 1}.

\textbf{Step 2.} \textit{Show that for all $p\geq2$, there exists a constant $C_{p,T}>0$ independent of $\Delta t$ and $N$ such that}
$\|X^N(T)-X^N_{n_{T}}\|_{L^p(\Omega,H)}\leq C_{p,T}(\Delta t)^{\frac{1-\delta}{2}}.$

On the partition $\mathcal{T}$, the process $\{X^N(t)\}_{t\in[0,T]}$ can be  rewritten as 
\begin{align}\label{eq:X^N_{N(T)}}
\begin{cases}
X^N(t_{i+1}-)\!=E(\Delta t_{i+1})X^N(t_{i})+\int_{t_i}^{t_{i+1}}\!\!E(t_{i+1}\!-\!s)F_N(X^N(s))ds+\int_{t_i}^{t_{i+1}}\!\!E(t_{i+1}\!-\!s)P_NdW(s)\\
\quad\quad \quad\quad\quad\quad -\int_{t_i}^{t_{i+1}}\int_\chi E(t_{i+1}\!-\!s)G_N(X^N(s),z)\nu(dz)ds,\\		X^N(t_{i+1})\!=X^N(t_{i+1}-)+G_N(X^N(t_{i+1}-),p_{t_{i+1}})\mathbf{1}_{\{p_{t_{i+1}}\neq 0\}}, \quad i=0,1,\ldots, n_{T}-1. 
\end{cases}
\end{align}
On the time interval $[t_i,t_{i+1})$, it follows from ~\eqref{fullnumericalscheme} and~\eqref{eq:X^N_{N(T)}} that
	 \begin{align*}
&\quad X^N_{i+1-}\!-\!X^N(t_{i+1}-)\nonumber\\
&=E(\Delta t_{i+1})(X^N_{i}\!-\!X^N(t_{i}))+\int_{t_{i}}^{t_{i+1}}E(t_{i+1}-s)\big(F_N(X^N_{i})\!-\!F_N(X^N(s))\big)ds\nonumber\\
&\quad-\int_{t_i}^{t_{i+1}}\!\!\int_\chi \big(E(t_{i+1}-t_i)G_N(X^N_{i},z)\!-\!E(t_{i+1}-s)G_N(X^N(s),z)\big)\nu(dz)ds
\end{align*}
and 
\[
X^N_{i+1}-X^N(t_{i+1})=X^N_{i+1-}-X^N(t_{i+1}-)+\big(G_N(X^N_{i+1-},p_{t_{i+1}})-G_N( X^N(t_{i+1}-),p_{t_{i+1}})\big)\mathbf{1}_{\{p_{t_{i+1}}\neq 0\}}. 
\]
Denote $e_{i}:=X^N_{i}-X^N(t_{i})$ and $e_{i-}:=X^N_{i-}-X^N(t_{i}-)$. Then it holds that 
\begin{align}\label{eq:e_i}
e_{i+1}=\big(\mathrm I+g_1(p_{t_{i+1}})\mathbf{1}_{\{p_{t_{i+1}}\neq 0\}}\big)e_{i+1-}=:\mathcal G_{i+1}e_{i+1-}.	
\end{align}
Using the definition of $G_N$ gives that for $i=0,1,\ldots,n_T-1,$
\begin{align}\label{eq:e_i+1}
e_{i+1-}=E(\Delta t_{i+1})\Gamma_i\mathcal{G}_ie_{i-}+\Xi(t_{i+1})+\sum_{j=1}^3\mathcal{R}_j(t_{i+1}),
\end{align} 
where $\Gamma_i:=\mathrm{I}-\Delta t_{i+1}\int_\chi g_1(z)\nu(dz)$, $
\Xi(t_{i+1}):=\int_{t_{i}}^{t_{i+1}}E(t_{i+1}\!-\!s)\big(F_N(X^N_{i})-F_N(X^N(t_{i}))\big)ds$, and 
\begin{align*}
	&\mathcal{R}_1(t_{i+1}):= \int_{t_{i}}^{t_{i+1}}E(t_{i+1}-s)\big(F_N(X^N(t_i))-F_N(X^N(s))\big)ds,\\
    &\mathcal R_2(t_{i+1}):=\int_{t_i}^{t_{i+1}}\int_\chi E(t_{i+1}-t_i)\big(G_N(X^N(s),z)-G_N(X^N(t_{i}),z)\big)\nu(dz)ds,\\
	&\mathcal R_3(t_{i+1}):=\int_{t_{i}}^{t_{i+1}}\int_\chi\big(E(t_{i+1}-s)-E(t_{i+1}-t_i)\big)G_N(X^N(s),z)\nu(dz)ds.
\end{align*}
By iteratively applying the relation between $e_{i+1-}$ and $e_{i-}$ (i.e., \eqref{eq:e_i+1}), the local error is accumulated into the global one, namely, we have
\begin{align*}
	e_{i+1-}=\sum_{k=0}^i\Big\{\prod_{l=k+1}^i E(\Delta t_{l+1})\Gamma_l\mathcal G_l\Big\}\Big(\Xi(t_{k+1})+\sum_{j=1}^3\mathcal{R}_j(t_{k+1})\Big),
\end{align*}
where we used $e_{0}=0$ and set $\prod_{i+1}^i:=1$. Taking $\|\cdot\|_{L^p(\Omega,H)}$-norm yields
\begin{align*}
	&\|e_{n_T-}\|_{L^p(\Omega,H)}
	\leq \Big\| \sum_{k=0}^{n_T-1}\Big\{\prod_{l=k+1}^{n_T-1}E(\Delta t_{l+1})\Gamma_l\mathcal{G}_l\Big\} \Xi(t_{k+1})\Big\|_{L^p(\Omega,H)}\\
	&+\Big\| \sum_{k=0}^{n_T-1}\Big\{\prod_{l=k+1}^{n_T-1}\!\!E(\Delta t_{l+1})\Gamma_l\mathcal{G}_l\Big\}\mathcal R_1(t_{k+1})\Big\|_{L^p(\Omega,H)}\!\!+\Big\| \sum_{k=0}^{n_T-1} \Big\{\prod_{l=k+1}^{n_T-1}\!\!E(\Delta t_{l+1})\Gamma_l\mathcal{G}_l\Big\}\mathcal R_2(t_{k+1})\Big\|_{L^p(\Omega,H)}\\
    &+\Big\|\sum_{k=0}^{n_T-1}\Big\{\prod_{l=k+1}^{n_T-1}E(\Delta t_{l+1})\Gamma_l\mathcal{G}_l\Big\}\mathcal R_3(t_{k+1})\Big\|_{L^p(\Omega,H)}=:\mathcal J_1+\mathcal J_2+\mathcal J_3+\mathcal J_4.
\end{align*} 
To proceed, we need to establish the $p$th moment estimates on terms $\mathcal J_1,\ldots,\mathcal J_4$ respectively, which involve the multiplicative noise or the jump process. 

Let $\Delta t_0=0$ and $\lfloor s\rfloor:=\min\{t_{n-1},n\in\{1,\ldots,n_T+1\}:\sum_{i=0}^{n}\Delta t_i>s~\mbox{or}~\sum_{i=0}^{n-1}\Delta t_i=s\}$ for any $s\in [0,T]$, and use $\ell(\lfloor s\rfloor)$ to denote the subscript corresponding to time $\lfloor s\rfloor$ (i.e., $t_{\ell(\lfloor s\rfloor)}=\lfloor s\rfloor$). For the term $\mathcal J_1$, we have
 \begin{align*}	
 \mathcal J_1&=\Big\|\int_0^{T}E(T-s)\Big\{\prod_{l=\ell(\lfloor s\rfloor)+1}^{n_T-1}\Gamma_l\mathcal{G}_l\Big\}\big(F_N(X^N_{\ell(\lfloor s\rfloor)})-F_N(X^N(\lfloor s\rfloor))\big)ds\Big\|_{L^p(\Omega,H)}\\
 &\leq \int_0^T\Big(\mathbb E\Big[\mathbb E\Big[\Big\|\prod_{l=\ell(\lfloor s\rfloor)+1}^{n_T-1}\Gamma_l\mathcal{G}_l\Big\|^p_{\mathcal L(H)}\big\|F_N(X^N_{\ell(\lfloor s\rfloor)})-F_N(X^N(\lfloor s\rfloor))\big\|^p\Big|\mathcal F_{\lfloor s\rfloor}\Big]\Big]\Big)^{\frac 1p}ds. 
 \end{align*}
Note that random variables $\Gamma_l$, $\mathcal G_l$, $l\ge \ell(\lfloor s\rfloor)+1$ and the $\sigma$-algebra $\mathcal F_{\lfloor s\rfloor}$  are independent, which together with the Lipschitz condition of $F$, \eqref{eq:e_i}, and $\|\mathcal G_l\|_{\mathcal L(H)}\leq C$  yields 
\begin{align*}
	\mathcal J_1 
 &\leq \int_0^T\Big\|\prod_{l=\ell(\lfloor s\rfloor)+1}^{n_T-1}\Gamma_l\mathcal{G}_l\Big\|_{L^p(\Omega,\mathcal L(H))}\big\|F_N(X^N_{\ell(\lfloor s\rfloor)})-F_N(X^N(\lfloor s\rfloor))\big\|_{L^p(\Omega,H)}ds\\
 &\leq L_F\Big(\E\Big[\sup_{k\in\{0,1,...,n_T-1\}}\prod_{l=k+1}^{n_T-1}\|\Gamma_l\mathcal{G}_l\|_{\IL(H)}^{p}\Big]\Big)^{\frac{1}{p}}\int_0^{T}\big\|\mathcal G_{\ell(\lfloor s\rfloor)}e_{\ell(\lfloor s\rfloor)-}\big\|_{L^{p}(\Omega,H)}ds. 
\end{align*}
According to Assumption~\ref{assum:jumpcoefficient} on coefficient $G$, it holds that $\|\Gamma_l\|_{\IL(H)}\leq 1+C_G\Delta t_{l+1}$ with $C_G:=b\nu(\chi)$, and $
\|\mathcal{G}_l\|_{\IL(H)}\leq 1+b\mathbf{1}_{\{p_{t_{l}}\neq 0\}}$. Set $N(T):=N((0,T]\times\chi)$. From Poisson distribution  
$\mathbb{P}(N(T)=m)=\frac{(T\nu(\chi))^m}{m!}e^{-T\nu(\chi)},$ it follows that for all $q\geq1$,
\begin{align}\label{ineq:Pi}
	&\E\Big[\sup_{k\in\{0,1,...,n_T-1\}}\prod_{l=k+1}^{n_T-1}\|\Gamma_l\mathcal{G}_l\|_{\IL(H)}^q\Big]\leq \E\Big[(1+b)^{qN(T)}\prod_{l=1}^{n_T}\big(1+C_G\Delta t_{l}\big)^{q}\Big]\nonumber\\
	&\leq \E\Big[(1+b)^{qN(T)}e^{qC_G\sum_{l=1}^{n_T}\Delta t_{l}}\Big]
	\leq C_{q,T}\sum_{m=1}^\infty (1+b)^{qm}\frac{(T\nu(\chi))^m}{m!}e^{-T\nu(\chi)}<\infty.
\end{align} 
Hence we arrive at $\mathcal J_1 \leq C_{p,T}\int_0^{T}\|e_{\ell(\lfloor s\rfloor)-}\|_{L^{p}(\Omega,H)}ds.$

For the term $\mathcal J_2,$ by Assumption \ref{assum:nonlinearity}, we have  
\begin{align*}
\mathcal J_2&= \Big\|\int_0^T\Big\{\prod_{l=\ell(\lfloor s\rfloor)+1}^{n_T-1}\!\! \Gamma_l\mathcal G_l\Big\}E(T\!-\!s)(-A)^{\frac{1+\delta}{4}}(-A)^{-\frac{1+\delta}{4}}\big(F_N(X^N(s))\!-\!F_N(X^N(\lfloor s\rfloor))\big)ds\Big\|_{L^p(\Omega,H)}\\
&\leq C\int_0^T\Big\|\prod_{l=\ell(\lfloor s\rfloor)+1}^{n_T-1}\!\!\Gamma_l\mathcal G_l\Big\|_{L^{2p}(\Omega,\mathcal L(H))}(T-s)^{-\frac{1+\delta}{4}}\big\|X^N(s)-X^N(\lfloor s\rfloor)\big\|_{L^{4p}(\Omega,\dot{H}^{-\frac{1-\delta}{2}})}\times \\
&\quad\quad\Big(1+\|X^N(s)\|_{L^{4p}(\Omega,\dot{H}^{\frac{1-\delta}{2}})}+\|X^N(\lfloor s\rfloor)\|_{L^{4p}(\Omega,\dot{H}^{\frac{1-\delta}{2}})}\Big)ds.  
\end{align*} 
To estimate $\|X^N(s)-X^N(\lfloor s\rfloor)\|_{L^{2q}(\Omega,\dot{H}^{-\frac{1-\delta}{2}})}$ for $q\geq1$ and $s\in[0,T]$,  we note that 
 \begin{align*}
 &\Big\|\int_{\lfloor s\rfloor}^{s}(-A)^{-\frac{1-\delta}{4}}E(s\!-\!r)P_NdW(r)\Big\|^{2q}_{L^{2q}(\Omega,H)}
  \leq \E\Big[\E\Big[\Big\|\int_{\zeta}^{s}\!(-A)^{-\frac{1-\delta}{4}}E(s\!-\!r)P_NdW(r)\Big\|^{2q}\Big]\Big|_{\zeta=\lfloor s\rfloor}\Big]\\
  &\leq C_q\E\Big[\Big(\int_{\zeta}^{s}\|(-A)^{\frac{\delta}{2}}(-A)^{-\frac{1+\delta}{4}}E(s-r)\|_{\IL_2(H)}^2dr\Big)^{q}\Big|_{\zeta=\lfloor s\rfloor}\Big]\\
  &\leq C_q\|(-A)^{-\frac{1+\delta}{4}}\|_{\IL_2(H)}^{2q}\E\big[(s-\lfloor s\rfloor)^{q(1-\delta)}\big]\leq C_q(\Delta t)^{q(1-\delta)}, 
\end{align*}
where we used the Burkholder--Davis--Gundy inequality (see e.g. \cite[Theorem 4.36]{GDP-2014}) and the H\"older inequality. This,  combining 
Proposition~\ref{prop:wellposedness} shows  
\begin{align*}
	&\quad\big\|X^N(s)-X^N(\lfloor s\rfloor)\big\|_{L^{2q}(\Omega,\dot{H}^{-\frac{1-\delta}{2}})}\\
	&\leq \Big\|(-A)^{-\frac{1-\delta}{4}}(E(s-\lfloor s\rfloor)-\mathrm{I})X^N(\lfloor s\rfloor)\Big\|_{L^{2q}(\Omega,H)}+\int_{(s-\Delta t)\vee 0}^{s}\big\|F_N(X^N(r))\big\|_{L^{2q}(\Omega,H)}dr\\
 &\quad+\Big\|\int_{\lfloor s\rfloor}^{s}(-A)^{-\frac{1-\delta}{4}}E(s-r)P_NdW(r)\Big\|_{L^{2q}(\Omega,H)}\\
 &\quad +\Big\|\int_{\lfloor s\rfloor}^{s}\int_\chi (-A)^{-\frac{1-\delta}{4}}E(s-r)G_N(X^N(r),z)\nu(dz)dr\Big\|_{L^{2q}(\Omega,H)}\\
	&\leq C_{q,T}(\Delta t)^{\frac{1-\delta}{2}}\!\!+\int_{(s-\Delta t)\vee 0}^s\int_\chi\Big(|g_1(z)|\big\|X^N(r)\big\|_{L^{2q}(\Omega,H)}+\|g(z)\|\Big)\nu(dz)dr\leq C_{q,T}(\Delta t)^{\frac{1-\delta}{2}}. 
\end{align*}
Hence, it follows from \eqref{ineq:Pi} that $\mathcal J_2\leq C_{p,T}(\Delta t)^{\frac{1-\delta}{2}}$.

For the term $\mathcal J_3,$ by \eqref{ineq:Pi}, similarly, we derive that
\begin{align*}
\mathcal J_3&\!\leq\! C\!\!\int_0^T\!\Big\|\prod_{l=\ell(\lfloor s\rfloor)+1}^{n_T-1} \!\!\!\!\!\Gamma_l\mathcal G_l\Big\|_{L^{2p}(\Omega,\mathcal L(H))}(T\!-\!\lfloor s\rfloor)^{-\frac{1-\delta}{4}}\big\|X^N(s)\!-\!X^N(\lfloor s\rfloor)\big\|_{L^{2p}(\Omega,\dot{H}^{-\frac{1-\delta}{2}})}ds\!\leq\! C_{p,T}(\Delta t)^{\frac{1-\delta}{2}}.
\end{align*}
The term $\mathcal J_4$ can be estimated as 
\begin{align*}
\mathcal J_4&\leq \int_0^T\!\!\!\int_{\chi}\Big\|\Big\{\prod_{l=\ell(\lfloor s\rfloor)+1}^{n_T-1}\!\!\!\!\!\Gamma_l\mathcal G_l\Big\}E(T\!-\!s)(-A)^{\frac12}(E(s\!-\!\lfloor s\rfloor)\!-\!\mathrm{I})(-A)^{-\frac12}G_N(X^N(s),z) \Big\|_{L^p(\Omega,H)}\nu(dz)ds \\&\leq  C(\Delta t)^{\frac12}\int_0^T\!\!\!\int_{\chi}(T-s)^{-\frac12} \big(|g_1(z)|\|X^N(s)\|_{L^{2p}(\Omega,H)}+\|g(z)\|\big)\nu(dz)ds\leq C_{p,T}(\Delta t)^{\frac12}. 
\end{align*}

Combining estimates of terms $\mathcal J_1,\ldots,\mathcal J_4$ above gives
\begin{align*}
\|e_{n_T-}\|_{L^p(\Omega,H)}\leq C_{p,T}\int_0^T\|e_{\ell(\lfloor s\rfloor)-}\|_{L^{p}(\Omega,H)}ds +C_{p,T}(\Delta t)^{\frac{1-\delta}{2}}, 
\end{align*}
which yields $\|e_{n_T-}\|_{L^p(\Omega,H)} \leq C_{p,T}(\Delta t)^{\frac{1-\delta}{2}}$ by using the Gr\"onwall inequality. According to \eqref{eq:e_i} and $\|\mathcal G_l\|_{\mathcal L(H)}\leq C$, we obtain
\begin{align*}
 \|e_{n_T}\|_{L^p(\Omega,H)}\leq C\|e_{n_T-}\|_{L^p(\Omega,H)}\leq C_{p,T}(\Delta t)^{\frac{1-\delta}{2}}.
\end{align*} 

Combining \textbf{Steps 1-2} finishes the proof. 
\end{proof}

\section{Proof of Theorem~\ref{theo:temporalstrong-Lp}.}\label{sec4}

In this section, we present the proof of the $L^p$-strong convergence orders of the fully discrete scheme \eqref{eq:fulldiscrete-onestep} for \eqref{SPDE:MultiplicativeNoise} driven by the additive Poisson noise with $\nu(\chi)\leq\infty$, that is Theorem \ref{theo:temporalstrong-Lp}. The proof relies on estimates of the stochastic convolutions, the recently obtained L\^e's quantitative John--Nirenberg inequality and some \textit{a priori} estimates on the nested conditional ``$L^2$-norms" of both stochastic convolutions and the solution of the perturbed stochastic equation.

 For $N\in\mathbb N_+$ and $t\in[0,T]$, denote stochastic convolutions
$\IW^N(t):=\int_0^{t}E(t-s)P_NdW(s)$ and  $
	\mathcal{N}^N(t):=\int_0^{t}\int_\chi E(t-s)g_N(z)\tilde{N}(dz,ds)$ with $g_N:=P_Ng$.
Let $Y^N(t):=X^N(t)-\IW^N(t)-\mathcal{N}^N(t)$, which satisfies the perturbed stochastic equation: 
\begin{align}\label{eq:perturbed}
dY^N(t)=AY^N(t)dt+F_N\big(Y^N(t)+\IW^{N}(t)+\mathcal N^N(t)\big)dt,\quad t\in(0,T],
\end{align}
with $Y^N(0)=P_Nx_0$. It is known that $Y^N$ satisfies
\begin{align}\label{purterb}
Y^N(t)=E(t)P_Nx_0+\int_0^tE(t-\rho)F_N(Y^N(\rho)+\IW^{N}(\rho)+\mathcal N^N(\rho))d\rho.
\end{align}
We list properties of the stochastic convolution with respect to the Wiener process as follows, whose proofs are postponed to Appendix.
\begin{lemma}
(\romannumeral1) \label{lem:WN-sup}
  For $k\geq 1$ and $0\leq s< t\leq T$, we have
  \begin{align}\label{ineq:WN-sup}
     \sup_{N\in\mathbb N_+}\mathbb E\Big[\sup_{t\in(s,T]}\Big\|\int_s^t(-A)^{\frac{1-\delta}{4}}E(t-r)P_NdW(r) &\Big\|^{2k}\Big]<\infty. 
  \end{align}
(\romannumeral2) 
For $p\geq1$ and $t\in[0,T]$, there exists a constant $C_{p,T}>0$ independent of $\Delta t$ such that
    \begin{align}\label{prop:localerror-WN}
       \sup_{N\in\mathbb N_+}\big\|(-A)^{-\frac{1-\delta}{4}}\big(\IW^N(t)-\IW^N(\lfloor t\rfloor)\big)\big\|_{L^{2p}(\Omega,H)}\leq C_{p,T}(\Delta t)^{\frac{1-\delta}{2}}. 
    \end{align}
\end{lemma}

We then introduce the recently obtained L\^e's quantitative John--Nirenberg inequality, which provides a useful way to 
 bound the $L^p$-norm of a stochastic process by moment estimates of the nested conditional ``$L^1$-norm". The conditional expectation given $\mathcal{F}_t$ is denoted by $\E^t$. 

\begin{lemma}\cite[Theorem 1.1]{KL2022}\label{ineq:WTNI}
  Let $(E,\mathrm d)$ be a metric space and $Z\colon[0,\infty)\times \Omega\to (E,\mathrm d)$ be a right continuous with left limits and adapted integrable stochastic process. Then for every $\tau>0$ and $p\geq1$, there exists a constant $C_p>0$ such that
  \begin{align*}
     \Big\|\sup_{t\in[0,\tau]}\mathrm d(Z_0,Z_t)\Big\|_{L^p(\Omega,\R)}\leq C_p\Big\|\sup_{t\in[0,\tau]}\E^t\big[\sup_{s\in[t,\tau]}\E^s[\mathrm d(Z_{s-},Z_{\tau})]\big]\Big\|_{L^p(\Omega,\R)}.
  \end{align*}
 \end{lemma}

 Let $[a,b]_{<}$ denote $\{(s,t)\in[a,b]^2\colon s<t\}$. For a random variable $\vartheta$, a sub-$\sigma$-algebra $\mathcal G\subset \mathcal F$, and $p\ge 1$, we set   
 $\|\vartheta\|_{L^p(\Omega,H)|\mathcal G}:=\big(\E[\|\vartheta\|^p\big|\mathcal G]\big)^{1/p}$. For a measurable mapping $f\colon [0,T]\times\Omega\to H$, we set $\|f\|_{\mathcal{C}_2^0|\mathcal{F}_s,[s,t]}:=\sup_{r\in[s,t]}\|f(r)\|_{L^2(\Omega,H)|\mathcal{F}_s}$ for any $s,t\in[0,T]_{<}$. Let $s':=\lfloor s\rfloor+\Delta t$, that is, $s'$ is the smallest grid point strictly bigger than $s$. 
 
The following propositions give some \textit{a priori} estimates including nested conditional ``$L^2$-norms" of both stochastic convolutions and the solution of the perturbed stochastic equation. The proofs are postponed to Section~\ref{sec5}.

\begin{proposition}\label{prop:NN}
For each $N\in\N_+$ and any time grid point $t_i$, $i=0,1,...,M$, we have
\begin{align}\label{prop:localerror-NN}
\big\|(-A)^{-\frac{1-\delta}{4}}\big(\mathcal{N}^N(\cdot)-\mathcal{N}^N(\lfloor \cdot\rfloor)\big)\big\|_{\mathcal{C}_2^0|\mathcal{F}_{t_i},[t_i,T]}
&\leq C(\Delta t)^{\frac{1-\delta}{2}}\Big(1+\|\mathcal N^N(t_i)\|_{\dot{H}^{\frac{1-\delta}{2}}}\Big).
\end{align}	
\end{proposition}

\begin{proposition}\label{prop:WNN-conditional}
For $p\ge 2$, there exists a constant $C_{p,T}>0$ such that
\begin{align}
   & \sup_{N\in\mathbb N_+}\big\|\sup_{t\in[0,T]}\E^t\big[\sup_{s\in[t,T]}\sup_{r\in[0,T]}\E^s\big[\|\mathcal N^N(r)\|\big]\big]\big\|_{L^p(\Omega,\R)}\leq C_{p,T},\label{ineq:23}\\
  &\sup_{N\in\mathbb N_+}\Big\|\sup_{t\in[0,T]}\E^t\Big[\sup_{s\in[t,T]}\E^s\big[\big\|(-A)^{\frac{1-\delta}{4}}\mathcal N^N(s')\big\|^2\big]\Big] \Big\|_{L^p(\Omega,\R)}\leq C_{p,T},\label{ineq:24}\\
  &\sup_{N\in\mathbb N_+}\Big\|\sup_{t\in[0,T]}\E^t\Big[\sup_{s\in[t,T]}\E^s\big[\big\|(-A)^{\frac{1-\delta}{4}}\mathcal N^N(\cdot)\big\|^2_{\mathcal{C}_2^0|\mathcal{F}_{s'},[0,T]}\big]\Big]\Big\|_{L^p(\Omega,\R)}\leq C_{p,T},\label{ineq:NNC}\\
 &\sup_{N\in\mathbb N_+}\Big\|\sup_{t\in[0,T]}\E^t\Big[\sup_{s\in[t,T]}\E^s\big[\big\|(-A)^{\frac{1-\delta}{4}}\IW^N(\cdot)\big\|^2_{\mathcal{C}_2^0|\mathcal{F}_{s'},[0,T]}\big]\Big]\Big\|_{L^p(\Omega,\R)}\leq C_{p,T}.\label{ineq:WNC}
  \end{align}  
\end{proposition}

\begin{proposition}\label{prop:YN}
(\romannumeral1)
  For $p\geq 2$ and $t\in[0,T]$, there exists a constant $C_{p,T}>0$ such that
\begin{align}\label{prop:localerror-YN}
 \sup_{N\in\mathbb N_+}\big\|(-A)^{-\frac{1-\delta}{4}}(Y^N(t)&-Y^N(\lfloor t\rfloor))\big\|_{L^{p}(\Omega,H)}\leq C_{p,T}(\Delta t)^{\frac{1-\delta}{2}}.
\end{align}

(\romannumeral2) 
For $p\geq2$, there exists a constant $C_{p,T}>0$ such that
\begin{align}\label{prop:YN-conditional}
      \sup_{N\in\mathbb N_+} \Big\|\sup_{t\in[0,T]}\E^t\Big[\sup_{s\in[t,T]}\E^s\big[\big\|(-A)^{\frac{1-\delta}{4}}Y^N(\cdot)\big\|^2_{\mathcal{C}_2^0|\mathcal{F}_{s'},[s',T]}\big]\Big]\Big\|_{L^p(\Omega,\R)}\leq C_{p,T}.	
\end{align}
\end{proposition}

With these preliminaries, we present the proof of Theorem~\ref{theo:temporalstrong-Lp}.

\begin{proof}[\textbf{Proof of Theorem~\ref{theo:temporalstrong-Lp}}]
The error between the exact solution and the semi-discrete numerical solution (i.e., $\|X(T)-X^N(T)\|_{L^p(\Omega,H)}\leq C_{p,T}N^{-\frac{1-\delta}{2}}$) can be proved by \textbf{Step 1} in the proof of Theorem \ref{the:convergenceorder-jumpadapted} and thus is omitted. Hence it suffices to estimate the error $\|X^N(T)-X^N_{n_T}\|_{L^p(\Omega,H)}$. 
 
 It follows from assumptions on $F$ and $g$, \eqref{prop:semigroup2}, and Lemma~\ref{ineq:maximalinequality} that 
\begin{align*}
&\|X^N(T)-X^N_{n_T}\|_{L^p(\Omega,H)}\leq \Big\|\int_0^{T}E(T-s)\big(F_N(X^N(s))-F_N(X^N_{\ell(\lfloor s\rfloor)})\big)ds\Big\|_{L^p(\Omega,H)}\\
	&\quad+\Big\|\int_0^{T}\int_\chi E(T-s)(\mathrm{I}-E(s-\lfloor s\rfloor))(-A)^{-\frac{1-\delta}{2}}(-A)^{\frac{1-\delta}{2}}g_N(z)\tilde{N}(dz,ds)\Big\|_{L^p(\Omega,H)}\\
	&\leq I_1+I_2+I_3+L_F\int_0^{T}\!\!\big\|X^N(\lfloor s\rfloor)\!-\!X^N_{\ell(\lfloor s\rfloor)}\big\|_{L^p(\Omega,H)}ds+C_{p,T}\Big(\int_0^{T}\!\!\big\|(\mathrm{I}\!-\!E(s\!-\!\lfloor s\rfloor))(-A)^{-\frac{1-\delta}{2}}\big\|^p_{\IL(H)}\\
&\quad\times\Big\{\int_\chi\big\|(-A)^{\frac{1-\delta}{2}}g_N(z)\big\|^p\nu(dz)\!+\!\Big(\int_\chi\big\|(-A)^{\frac{1-\delta}{2}}g_N(z)\big\|^2\nu(dz)\Big)^{\frac{p}{2}}\Big\}ds\Big)^{\frac{1}{p}}\\
&\leq C_{p,T}(\Delta t)^{\frac{1-\delta}{2}}+L_F\int_0^{T}\big\|X^N(\lfloor s\rfloor)\!-\!
X^N_{\ell(\lfloor s\rfloor)}\big\|_{L^p(\Omega,H)}ds+I_1+I_2+I_3,
\end{align*}
where 
\begin{align*}
&I_1\!:=\!\!\Big\|\int_0^{T}E(T\!-\!s)\big(F_N(Y^N(s)+\IW^N(\lfloor s\rfloor)+\mathcal{N}^N(\lfloor s\rfloor))-F_N(X^N(\lfloor s\rfloor)\big)\big)ds\Big\|_{L^p(\Omega,H)},\\
&I_2\!:=\!\!\Big\|\int_0^{T}\!\!\!E(T\!-\!s)\big(F_N(Y^N(s)\!+\!\IW^N(s)\!+\!\mathcal{N}^N(\lfloor s\rfloor))\!-\!F_N(Y^N(s)\!+\!\IW^N(\lfloor s\rfloor)\!+\!\mathcal{N}^N(\lfloor s\rfloor))\big)ds\Big\|_{L^p(\Omega,H)},\\
  &I_3\!:=\!\!\Big\|\int_0^{T}E(T-s)\big(F_N(X^N(s))-F_N(Y^N(s)+\IW^N(s)+\mathcal{N}^N(\lfloor s\rfloor))\big)ds\Big\|_{L^p(\Omega,H)}.\end{align*} 
  
For the term $I_1$, by Assumption~\ref{assum:nonlinearity}, we have
\begin{align*}
  &I_1\leq C\int_0^{T}\big\|E(T-s)(-A)^{\frac{1+\delta}{4}}\big\|_{\IL(H)}\big\|(-A)^{-\frac{1-\delta}{4}}(Y^N(s)-Y^N(\lfloor s\rfloor))\big\|_{L^{2p}(\Omega,H)}ds~ \times \\
  &~~~ \Big(1+\sup_{s\in[0,T]}\|Y^N(s)\|_{L^{2p}(\Omega,\dot{H}^{\frac{1-\delta}{2}})}+\sup_{s\in[0,T]}\|\IW^N(s)\|_{L^{2p}(\Omega,\dot{H}^{\frac{1-\delta}{2}})}+\sup_{s\in[0,T]}\|\mathcal N^N(s)\|_{L^{2p}(\Omega,\dot{H}^{\frac{1-\delta}{2}})}\Big).
  \end{align*}
  Combining \eqref{prop:localerror-YN} and the same arguments as in the proof of Proposition~\ref{prop:wellposedness} leads to
\begin{align*}
  I_1\leq C_{p,T}(\Delta t)^{\frac{1-\delta}{2}}\int_0^{T}(T-s)^{-\frac{1+\delta}{4}}ds\leq C_{p,T}(\Delta t)^{\frac{1-\delta}{2}}.
\end{align*} 
Similarly, by \eqref{prop:localerror-WN},  we have
\begin{align*}
 I_2 &\leq C\int_0^{T}\!\!(T\!-\!s)^{-\frac{1+\delta}{4}}\big\|(-A)^{-\frac{1-\delta}{4}}\big(\IW^N(s)\!-\!\IW^N(\lfloor s\rfloor)\big)\big\|_{L^{2p}(\Omega,H)}ds\times\Big(1\!+\!\!\sup_{s\in[0,T]}\!\!\|Y^N(s)\|_{L^{2p}(\Omega,\dot{H}^{\frac{1-\delta}{2}})}\!\\
 &\quad +\!\sup_{s\in[0,T]}\|\IW^N(s)\|_{L^{2p}(\Omega,\dot{H}^{\frac{1-\delta}{2}})}\!+\!\sup_{s\in[0,T]}\|\mathcal N^N(s)\|_{L^{2p}(\Omega,\dot{H}^{\frac{1-\delta}{2}})}\Big)\leq  C_{p,T}(\Delta t)^{\frac{1-\delta}{2}}.
\end{align*}

For the term $I_3$, we claim that  
\begin{align}\label{estimate:I3}
	I_3\leq C_{p,T}(\Delta t)^{\frac{1-\delta}{2}}.
\end{align}	
Estimate of term $I_3$ includes dealing with the temporal H\"older continuity of the stochastic convolution with respect to the compensated Poisson random measure to overcome the order barrier. To this end we aim to apply Lemma~\ref{ineq:WTNI} to estimate the $L^p$-norm of $I_3$.  For $(s,t)\in[0,T]_{<}$, define an integrable  $\mathcal{F}_t$-adapted stochastic process:
\begin{align*}
\mathcal{A}^{N,\Delta}(t)&:=\int_0^{t}E(T-r)\big(F_N(X^N(r))-F_N(Y^N(r)+\IW^N(r)+\mathcal{N}^N(\lfloor r\rfloor))\big)dr 
\end{align*}
with $\mathcal{A}^{N,\Delta}(0)=0$. By the globally Lipschitz condition of $F$, we have
\begin{align*}
	\E\big[\big\|\mathcal{A}^{N,\Delta}(t)-\mathcal{A}^{N,\Delta}(s)\big\|^2\big]
	&\leq L_F|t-s|\int_s^{t}\E\big[\big\|\mathcal{N}^N(r)-\mathcal{N}^N(\lfloor r\rfloor)\big\|^2\big]dr\\
&\leq C|t-s|^2\sup_{r\in[0,T]}\E[\|\mathcal{N}^N(r)\|^2]\leq C_T|t-s|^2. 
\end{align*} 
The Kolmogorov continuity theorem implies that $\mathcal{A}^{N,\Delta}$ is a.s. continuous. 
Thus by Lemma~\ref{ineq:WTNI}, we arrive at  
\begin{align*}
  \big\|\sup_{t\in[0,T]}\|\mathcal{A}^{N,\Delta}(t)\|\big\|_{L^p(\Omega,\R)}\leq C_p\Big\|\sup_{t\in[0,T]}\E^t\Big[\sup_{s\in[t,T]}\E^s\big[\big\|\mathcal{A}^{\mathcal{N},\Delta}(T)-\mathcal{A}^{\mathcal{N},\Delta}(s)\big\|\big]\Big]\Big\|_{L^p(\Omega,\R)},
\end{align*}
which yields that
\begin{align}\label{ineq:Lp-convergence}
   I_3=\|\mathcal A^{N,\Delta}(T)\|_{L^p(\Omega,H)}\leq C_p\Big\|\sup_{t\in[0,T]}\E^t\Big[\sup_{s\in[t,T]}\E^s\big[\big\|\mathcal{A}^{\mathcal{N},\Delta}(T)-\mathcal{A}^{\mathcal{N},\Delta}(s)\big\|\big]\Big]\Big\|_{L^p(\Omega,\R)}.
\end{align}
Now we turn to estimating the term on the right hand side of the inequality in \eqref{ineq:Lp-convergence}. We first show that, there exists a real-valued adapted stochastic process $\xi$  such that
	\begin{align}\label{findxi}\E^s\big[\big\|\mathcal{A}^{N,\Delta}(T)-\mathcal{A}^{N,\Delta}(s)\big\|\big]\leq \xi_s\quad a.s.
	\end{align}  
Indeed, using Assumption~\ref{assum:nonlinearity} and   the conditional H\"older inequality gives
\begin{align*}
	& \E^s\big[\big\|\mathcal{A}^{N,\Delta}(T)-\mathcal{A}^{N,\Delta}(s)\big\|\big]\leq C\int_s^T(T-r)^{-\frac{1+\delta}{4}}\big\|(-A)^{-\frac{1-\delta}{4}}\big(\mathcal{N}^N(r)-\mathcal{N}^N(\lfloor  r\rfloor)\big)\big\|_{L^2(\Omega,H)|\mathcal{F}_s}dr\\
	&\times \Big(1\!+\!\big\|(-A)^{\frac{1-\delta}{4}}Y^N(\cdot)\big\|_{\mathcal{C}_2^0|\mathcal{F}_s,[s,T]}\!+\!\big\|(-A)^{\frac{1-\delta}{4}}\IW^N(\cdot)\big\|_{\mathcal{C}_2^0|\mathcal{F}_s,[s,T]}\!+\!\big\|(-A)^{\frac{1-\delta}{4}}\mathcal N^N(\cdot)\big\|_{\mathcal{C}_2^0|\mathcal{F}_s,[0,T]}\Big)\\
	&\leq C_T\big\|(-A)^{-\frac{1-\delta}{4}}\big(\mathcal{N}^N(\cdot)-\mathcal{N}^N(\lfloor \cdot\rfloor)\big)\big\|_{\mathcal{C}_2^0|\mathcal{F}_s,[s,T]}\times\\
    &\quad \Big(1\!+\!\big\|(-A)^{\frac{1-\delta}{4}}Y^N(\cdot)\big\|_{\mathcal{C}_2^0|\mathcal{F}_s,[s,T]}\!+\!\big\|(-A)^{\frac{1-\delta}{4}}\IW^N(\cdot)\big\|_{\mathcal{C}_2^0|\mathcal{F}_s,[s,T]}\!+\!\big\|(-A)^{\frac{1-\delta}{4}}\mathcal N^N(\cdot)\big\|_{\mathcal{C}_2^0|\mathcal{F}_s,[0,T]}\Big).
\end{align*}	
 In the case of $s\in\{0,\Delta t, 2\Delta t,...,T\}$, i.e., $\lfloor s\rfloor=s$, it follows from \eqref{prop:localerror-NN} and the Young inequality that $\E^s\big[\|\mathcal{A}^{N,\Delta}(T)-\mathcal{A}^{N,\Delta}(s)\|\big]\leq \xi^1_s$ almost surely with 
\begin{align*}
   \xi^1_s:&= C_T(\Delta t)^{\frac{1-\delta}{2}}\Big(1+\big\|(-A)^{\frac{1-\delta}{4}}\mathcal N^N(s)\big\|^2+\big\|(-A)^{\frac{1-\delta}{4}}Y^N(\cdot)\big\|^2_{\mathcal{C}_2^0|\mathcal{F}_s,[s,T]}\\
   &\quad +\big\|(-A)^{\frac{1-\delta}{4}}\IW^N(\cdot)\big\|^2_{\mathcal{C}_2^0|\mathcal{F}_s,[s,T]}+\big\|(-A)^{\frac{1-\delta}{4}}\mathcal N^N(\cdot)\big\|^2_{\mathcal{C}_2^0|\mathcal{F}_s,[0,T]}\Big). 
\end{align*}
In the case that $s$ is not a grid point, i.e., $\lfloor s\rfloor\neq s$, we note that when $T-s\leq \Delta t$, 
\begin{align*}
 \E^s\big[\big\|\mathcal{A}^{N,\Delta}(T)\!-\!\mathcal{A}^{N,\Delta}(s)\big\|\big]
  &\leq C\Delta t\!\!\sup_{r\in[s,T]}\!\E^s\big[\big\|\mathcal N^N(r)\!-\!\mathcal{N}^N(\lfloor  r\rfloor)\big\|\big]\leq C\Delta t\!\!\sup_{r\in[s\!-\!\Delta t,T]}\!\!\E^s\big[\big\|\mathcal N^N(r)\big\|\big].
\end{align*}
When $T-s> \Delta t$, we observe that $s'-s\leq \Delta t$, $s'\leq T$, and
\begin{align*}
  \E^s\big[\big\|\mathcal{A}^{N,\Delta}(T)-\mathcal{A}^{N,\Delta}(s)\big\|\big]   &\leq \E^s\big[\E^{s'}\big[\|\mathcal{A}^{N,\Delta}(T)-\mathcal{A}^{N,\Delta}(s')\|\big]\big]+ \E^s\big[\|\mathcal{A}^{N,\Delta}(s)-\mathcal{A}^{N,\Delta}(s')\|\big]\\
 &\leq \E^s\big[\xi^1_{s'}\big]+\E^s\big[\|\mathcal{A}^{N,\Delta}(s)-\mathcal{A}^{N,\Delta}(s')\|\big].
\end{align*} 
Hence, we derive that $\E^s\big[\|\mathcal{A}^{N,\Delta}(T)-\mathcal{A}^{N,\Delta}(s)\|\big]\leq \xi^2_s$ almost surely with  
\begin{align*}
 \xi_s^2&:=C_T(\Delta t)^{\frac{1-\delta}{2}}\Big(1+\E^s\Big[\big\|(-A)^{\frac{1-\delta}{4}}\mathcal N^N(s')\big\|^2+\big\|(-A)^{\frac{1-\delta}{4}}Y^N(\cdot)\big\|^2_{\mathcal{C}_2^0|\mathcal{F}_{s'},[s',T]}\\
   &\quad +\big\|(-A)^{\frac{1-\delta}{4}}\IW^N(\cdot)\big\|^2_{\mathcal{C}_2^0|\mathcal{F}_{s'},[s',T]}\!+\!\big\|(-A)^{\frac{1-\delta}{4}}\mathcal N^N(\cdot)\big\|^2_{\mathcal{C}_2^0|\mathcal{F}_{s'},[0,T]}\Big]\!+\!\!\sup_{r\in[s-\Delta t,T]}\!\E^s\big[\|\mathcal N^N(r)\|\big]\Big). 
\end{align*}
Therefore \eqref{findxi} holds with $\xi_s:=\xi_s^1\mathbf{1}_{\{s=\lfloor s\rfloor\}}+\xi_s^2\mathbf{1}_{\{s\neq\lfloor s\rfloor\}}$. To further estimate the term on the right hand side of \eqref{ineq:Lp-convergence}, we use Propositions \ref{prop:NN}--\ref{prop:YN} to obtain
\begin{align*}
 & \big\|\sup_{t\in[0,T]}\E^t\big[\sup_{s\in[t,T]}\xi^2_s\big] \big\|_{L^p(\Omega,\R)}
 \leq C_T(\Delta t)^{\frac{1-\delta}{2}}\Big\{1+\Big\|\sup_{t\in[0,T]}\E^t\Big[\sup_{s\in[t,T]}\E^s\big[\big\|(-A)^{\frac{1-\delta}{4}}\mathcal N^N(s')\big\|^2\big]\Big]\Big\|_{L^p(\Omega,\R)}\\
 &~~~~~~~~~~~~~~~+\big\|\sup_{t\in[0,T]}\E^t\big[\sup_{s\in[t,T]}\sup_{r\in[s-\Delta t,T]}\E^s\big[\|\mathcal N^N(r)\|\big]\big]\big\|_{L^p(\Omega,\R)}\\
 &~~~~~~~~~~~~~~~+\Big\|\sup_{t\in[0,T]}\E^t\Big[\sup_{s\in[t,T]}\E^s\big[\big\|(-A)^{\frac{1-\delta}{4}}\mathcal N^N(\cdot)\big\|^2_{\mathcal{C}_2^0|\mathcal{F}_{s'},[0,T]}\big]\Big]\Big\|_{L^p(\Omega,\R)}\\
&~~~~~~~~~~~~~~~+\Big\|\sup_{t\in[0,T]}\E^t\Big[\sup_{s\in[t,T]}\E^s\big[\big\|(-A)^{\frac{1-\delta}{4}}\IW^N(\cdot)\big\|^2_{\mathcal{C}_2^0|\mathcal{F}_{s'},[s',T]}\big]\Big]\Big\|_{L^p(\Omega,\R)}\\
 &~~~~~~~~~~~~~~~+\Big\|\sup_{t\in[0,T]}\E^t\Big[\sup_{s\in[t,T]}\E^s\big[\big\|(-A)^{\frac{1-\delta}{4}}Y^N(\cdot)\big\|^2_{\mathcal{C}_2^0|\mathcal{F}_{s'},[s',T]}\big]\Big]\Big\|_{L^p(\Omega,\R)}\Big\}\leq C_{p,T}(\Delta t)^{\frac{1-\delta}{2}}.
\end{align*}
Similarly, we can derive that
$\big\|\sup_{t\in[0,T]}\E^t\big[\sup_{s\in[t,T]}\xi^1_s\big] \big\|_{L^p(\Omega,\R)}\leq C_{p,T}(\Delta t)^{\frac{1-\delta}{2}}.$
As a result, \eqref{estimate:I3} is proved. 

Combining estimates of terms $I_1$, $I_2$, and $I_3$, and applying the Gr\"onwall inequality complete the proof of the theorem.
\end{proof}

\begin{remark}
    In Theorem~\ref{theo:temporalstrong-Lp}, if we get rid of condition~\eqref{condition:g} imposed on the jump coefficient $g$, then we obtain that 
  $
\|X^N(T)-X^N_M\|_{L^p(\Omega,H)}\leq C_{p,T}(\Delta t)^{\frac{1}{4}},\; p\ge 2 
$
under Assumptions~\ref{assum:nonlinearity}--\ref{assum:x0} and $g_1\equiv0$. This is due to that the estimate of $\|\int_0^{T}\!\!\int_\chi (E(T\!-s)-E(T\!-\!\lfloor s\rfloor))g_N(z)\tilde{N}(dz,ds)\|_{L^p(\Omega,H)}$ can only provide the $L^p$-strong  convergence order $\frac14$ for all $p\geq2$.
\end{remark}

\section{Proofs of some essential propositions}\label{sec5}
In this section, we give proofs of Propositions \ref{prop:NN}--\ref{prop:YN}. Note that the process $\int^{\cdot}_t\int_\chi E(\cdot-r)g_N(z)\tilde{N}(dz,dr)$ is progressively measurable on product space $(\Omega\times[t,T], \mathcal F_T\times\mathcal B([t,T]), \mathbb P\times dr)$.
 For any interval $[a,b]\in\mathcal B([t,T])$, due to the independent increments property of the compensated Poisson measure, it holds that 
\begin{itemize}
\item[(P1)] $\sup\limits_{s\in [a,b]}\Big\|\int^s_t\int_\chi E(s-r)(-A)^{\frac{1-\delta}{4}}g_N(z)\tilde{N}(dz,dr)\Big\| $ is $\mathcal F_t$-independent.
\end{itemize}
This property will be frequently utilized in the following proofs.

\begin{proof}[\textbf{Proof of Proposition \ref{prop:NN}}]
For any time grid point $t_i$, $i=0,1,\ldots,M$, the property (P1) and the assumption on $g$ lead to
\begin{align*}
&\quad\big\|(-A)^{-\frac{1-\delta}{4}}\big(\mathcal{N}^N(\cdot)-\mathcal{N}^N(\lfloor \cdot\rfloor)\big)\big\|_{\mathcal{C}_2^0|\mathcal{F}_{t_i},[t_i,T]}\\
&\leq \sup_{s\in[t_i,T]} \Big\|\int_{\lfloor  s\rfloor}^s\int_\chi (-A)^{-\frac{1-\delta}{4}}E(s-\rho)g_N(z)\tilde{N}(dz,d\rho)\Big\|_{L^2(\Omega,H)}\\
 &\quad+\sup_{s\in[t_i,T]}\Big\|\int_0^{\lfloor  s\rfloor}\!\!\!\int_\chi (-A)^{-\frac{1-\delta}{4}}\big(E(s\!-\!\rho)-E(\lfloor s\rfloor\!-\!\rho)\big)g_N(z)\tilde{N}(dz,d\rho)\Big\|_{L^2(\Omega,H)|\mathcal{F}_{t_i}}\\
&\leq C(\Delta t)^{\frac12}\!+\!\!\sup_{s\in[t_i,T]}\!\Big\|\big(E(s\!-\!\lfloor s\rfloor)\!-\!\mathrm{I}\big)(-A)^{-\frac{1-\delta}{2}}\Big\|_{\IL(H)}\Big\|\!\!\int^{t_i}_0\!\!\!\int_\chi (-A)^{\frac{1-\delta}{4}}E(\lfloor s\rfloor\!-\!\rho)g_N(z)\tilde{N}(dz,d\rho)\Big\|\\
 &\quad+\!\!\sup_{s\in[t_i,T]}\Big\|\big(E(s\!-\!\lfloor s\rfloor)\!-\mathrm{I}\big)(-A)^{-\frac{1-\delta}{2}}\Big\|_{\IL(H)}\Big\|\!\!\int_{t_i}^{\lfloor  s\rfloor}\!\!\!\int_\chi(-A)^{\frac{1-\delta}{4}}E(\lfloor  s\rfloor-\rho)g_N(z)\tilde{N}(dz,d\rho)\Big\|_{L^2(\Omega,H)}\\
 &\leq C(\Delta t)^{\frac{1-\delta}{2}}\Big(1+\Big\|\int^{t_i}_0\int_\chi E(t_i-s)g_N(z)\tilde{N}(dz,ds)\Big\|_{\dot{H}^{\frac{1-\delta}{2}}}\Big).
\end{align*}	
The proof is completed.
\end{proof}

\begin{proof}[\textbf{Proof of Proposition~
\ref{prop:WNN-conditional}}]
(i)  \textit{The proof of \eqref{ineq:23}.}
It follows from the property (P1) that 
\begin{align*}
  &\quad\big\|\sup_{t\in[0,T]}\E^t\big[\sup_{s\in[t,T]}\sup_{r\in[0,T]}\E^s\big[\|\mathcal N^N(r)\|\big]\big]\big\|_{L^p(\Omega,\R)}\\
  &\leq \Big\|\sup_{t\in[0,T]}\E^t\Big[\sup_{s\in[t,T]}\sup_{r\in[s,T]}\E^s\Big[\Big\|\Big(\int^s_0+\int_{s}^{r}\Big)\int_\chi E(r-\rho)g_N(z)\tilde{N}(dz,d\rho)\Big\|\Big]\Big]\Big\|_{L^p(\Omega,\R)}\\
  &\quad+\Big\|\sup_{t\in[0,T]}\E^t\Big[\sup_{s\in[t,T]}\sup_{r\in[0,s]}\E^s\Big[\Big\|\int^{r}_0\int_\chi E(r-\rho)g_N(z)\tilde{N}(dz,d\rho)\Big\|\Big]\Big]\Big\|_{L^p(\Omega,\R)} \\
  &\leq \Big\|\sup_{t\in[0,T]}\E^t\Big[\sup_{s\in[t,T]}\sup_{r\in[s,T]}\E\Big[\Big\|\int^r_s\int_\chi E(r-\rho)g_N(z)\tilde{N}(dz,d\rho)\Big\|\Big]\Big]\Big\|_{L^p(\Omega,\R)}\\
  &\quad+\Big\|\sup_{t\in[0,T]}\E^t\Big[\sup_{s\in[t,T]}\Big\|\Big(\int^t_0+\int_t^s\Big)\int_\chi E(s-\rho)g_N(z)\tilde{N}(dz,d\rho)\Big\|\Big]\Big\|_{L^p(\Omega,\R)}\\
  &\quad+\Big\|\sup_{t\in[0,T]}\E^t\Big[\sup_{r\in[0,T]}\Big\|\int^{r}_0\int_\chi E(r-\rho)g_N(z)\tilde{N}(dz,d\rho)\Big\|\Big]\Big\|_{L^p(\Omega,\R)},
  \end{align*}
 and further by  Lemma~\ref{ineq:maximalinequality}, we have 
  \begin{align*}
&\quad\big\|\sup_{t\in[0,T]}\E^t\big[\sup_{s\in[t,T]}\sup_{r\in[0,T]}\E^s\big[\|\mathcal N^N(r)\|\big]\big]\big\|_{L^p(\Omega,\R)}\\
  &\leq C_{p,T}+\Big\|\sup_{t\in[0,T]}\Big\|\int^t_0\int_\chi E(t-\rho)g_N(z)\tilde{N}(dz,d\rho)\Big\|\Big\|_{L^p(\Omega,\R)}\\
  &\quad+\Big\|\sup_{t\in[0,T]}\E\Big[\sup_{s\in[t,T]}\Big\|\int_t^s\int_\chi E(s-\rho)g_N(z)\tilde{N}(dz,d\rho)\Big\|\Big]\Big\|_{L^p(\Omega,\R)}\\
  &\quad+\Big\|\sup_{t\in[0,T]}\E^t\Big[\sup_{r\in[0,t]}\Big\|\int_0^{r}\int_\chi E(r-\rho)g_N(z)\tilde{N}(dz,d\rho)\Big\|\Big]\Big\|_{L^p(\Omega,\R)} \\
  &\quad+\Big\|\sup_{t\in[0,T]}\E^t\Big[\sup_{r\in[t,T]}\Big\|\Big(\int^{t}_0+\int^{r}_t\Big)\int_\chi E(r-\rho)g_N(z)\tilde{N}(dz,d\rho)\Big\|\Big]\Big\|_{L^p(\Omega,\R)}\\
&\leq C_{p,T}+\Big\|\sup_{r\in[0,T]}\Big\|\int^{r}_0\int_\chi E(r-\rho)g_N(z)\tilde{N}(dz,d\rho)\Big\|\Big\|_{L^p(\Omega,\R)}\\
  &\quad+\Big\|\sup_{t\in[0,T]}\E\Big[\sup_{r\in[t,T]}\Big\|\int^{r}_t\int_\chi E(r-\rho)g_N(z)\tilde{N}(dz,d\rho)\Big\|\Big]\Big\|_{L^p(\Omega,\R)}\leq C_{p,T}.
\end{align*}

(ii) \textit{The proof of \eqref{ineq:24}.}
The proof is similar to that of \eqref{ineq:23}. We have
\begin{align*}
  &\quad\Big\|\sup_{t\in[0,T]}\E^t\Big[\sup_{s\in[t,T]}\E^s\big[\big\|(-A)^{\frac{1-\delta}{4}}\mathcal N^N(s')\big\|^2\big]\Big] \Big\|_{L^p(\Omega,\R)}\\
   &\leq C\Big\|\sup_{t\in[0,T]}\E^t\Big[\sup_{s\in[t,T]}\E\Big[\Big\|\int^{s'}_s\int_\chi E(s'-r)(-A)^{\frac{1-\delta}{4}}g_N(z)\tilde{N}(dz,dr)\Big\|^2\Big] \Big\|_{L^p(\Omega,\R)}\\
   &\quad+C\Big\|\sup_{t\in[0,T]}\E^t\Big[\sup_{s\in[t,T]}\Big\|\int^s_0\int_\chi E(s'-r)(-A)^{\frac{1-\delta}{4}}g_N(z)\tilde{N}(dz,dr)\Big\|^2\Big]\Big\|_{L^p(\Omega,\R)}\\
   &\leq C\Big\|\sup_{t\in[0,T]}\E^t\Big[\sup_{s\in[t,T]}\int^{s'}_s\int_\chi \big\|E(s'-r)(-A)^{\frac{1-\delta}{4}}g_N(z)\big\|^2\nu(dz)dr\Big] \Big\|_{L^p(\Omega,\R)}\\
   &\quad+C\Big\|\sup_{t\in[0,T]}\E^t\Big[\sup_{s\in[t,T]}\Big\|\Big(\int^s_t+\int^t_0\Big)\int_\chi E(s'-r)(-A)^{\frac{1-\delta}{4}}g_N(z)\tilde{N}(dz,dr)\Big\|^2\Big]\Big\|_{L^p(\Omega,\R)}\\
   &\leq C_{T}+C\Big\|\sup_{t\in[0,T]}\E\Big[\sup_{s\in[t,T]}\Big\|\int^s_t\int_\chi E(s-r)(-A)^{\frac{1-\delta}{4}}g_N(z)\tilde{N}(dz,dr)\Big\|^2\Big]\Big\|_{L^p(\Omega,\R)}\\
   &\quad+C\Big\|\sup_{t\in[0,T]}\Big\|\int^t_0\int_\chi E(t-r)(-A)^{\frac{1-\delta}{4}}g_N(z)\tilde{N}(dz,dr)\Big\|^2\Big\|_{L^p(\Omega,\R)}\leq C_{p,T}.
\end{align*}

\begin{itemize}
	\item[(iii)] \textit{The proof of \eqref{ineq:NNC}.} Based on the property (P1), we have
\end{itemize}
\begin{align*}
&\quad\Big\|\sup_{t\in[0,T]}\E^t\Big[\sup_{s\in[t,T]}\E^s\big[\big\|(-A)^{\frac{1-\delta}{4}}\mathcal N^N(\cdot)\big\|^2_{\mathcal{C}_2^0|\mathcal{F}_{s'},[0,T]}\big]\Big]\Big\|_{L^p(\Omega,\R)}\\
  &\leq C\Big\|\sup_{t\in[0,T]}\E^t\Big[\sup_{s\in[t,T]}\E^s\Big[\sup_{r\in[0,s']}\Big\|\int^r_0\int_\chi E(r-\rho)(-A)^{\frac{1-\delta}{4}}g_N(z)\tilde{N}(dz,d\rho)\Big\|^2\Big]\Big]\Big\|_{L^p(\Omega,\R)} \\
  &+\!C\Big\|\sup_{t\in[0,T]}\E^t\Big[\sup_{s\in[t,T]}\E^s\Big[\sup_{r\in[s',T]}\E^{s'}\Big[\Big\|\Big(\int^r_{s'}\!+\!\int_0^{s'}\Big)\!\int_\chi E(r\!-\!\rho)(-A)^{\frac{1-\delta}{4}}g_N(z)\tilde{N}(dz,d\rho)\Big\|^2\Big]\Big]\Big]\Big\|_{L^p(\Omega,\R)}\\
  &\leq C\Big\|\sup_{t\in[0,T]}\E^t\Big[\sup_{s\in[t,T]}\E^s\Big[\sup_{r\in[0,s]}\Big\|\int^r_0\int_\chi E(r\!-\!\rho)(-A)^{\frac{1-\delta}{4}}g_N(z)\tilde{N}(dz,d\rho)\Big\|^2\Big]\Big]\Big\|_{L^p(\Omega,\R)} \\
  &\quad+C\Big\|\sup_{t\in[0,T]}\E^t\Big[\sup_{s\in[t,T]}\E^s\Big[\sup_{r\in[s,s']}\Big\|\Big(\int_0^s+\int^r_s\Big)\int_\chi E(r\!-\!\rho)(-A)^{\frac{1-\delta}{4}}g_N(z)\tilde{N}(dz,d\rho)\Big\|^2\Big]\Big]\Big\|_{L^p(\Omega,\R)} \\
  &\quad+C\Big\|\sup_{t\in[0,T]}\E^t\Big[\sup_{s\in[t,T]}\E^s\Big[\sup_{r\in[s',T]}\E\Big[\Big\|\int^r_{s'}\int_\chi E(r\!-\!\rho)(-A)^{\frac{1-\delta}{4}}g_N(z)\tilde{N}(dz,d\rho)\Big\|^2\Big]\Big]\Big]\Big\|_{L^p(\Omega,\R)}\\
  &\quad+C\Big\|\sup_{t\in[0,T]}\E^t\Big[\sup_{s\in[t,T]}\E^s\Big[\Big\|\Big(\int_0^s+\int_s^{s'}\Big)\int_\chi E(s'\!-\!\rho)(-A)^{\frac{1-\delta}{4}}g_N(z)\tilde{N}(dz,d\rho)\Big\|^2\Big]\Big]\Big\|_{L^p(\Omega,\R)}.
  \end{align*}
  And using Lemma~\ref{ineq:maximalinequality} again, we obtain
  \begin{align*}
   &\quad\Big\|\sup_{t\in[0,T]}\E^t\Big[\sup_{s\in[t,T]}\E^s\big[\big\|(-A)^{\frac{1-\delta}{4}}\mathcal N^N(\cdot)\big\|^2_{\mathcal{C}_2^0|\mathcal{F}_{s'},[0,T]}\big]\Big]\Big\|_{L^p(\Omega,\R)} \\
   &\leq C_{p,T}+C\Big\|\sup_{t\in[0,T]}\E^t\Big[\sup_{r\in[0,t]}\Big\|\int^r_0\int_\chi E(r-\rho)(-A)^{\frac{1-\delta}{4}}g_N(z)\tilde{N}(dz,d\rho)\Big\|^2\Big]\Big\|_{L^p(\Omega,\R)} \\
  &\quad+C\Big\|\sup_{t\in[0,T]}\E^t\Big[\sup_{r\in[t,T]}\Big\|\Big(\int_0^t+\int^r_t\Big)\int_\chi E(r-\rho)(-A)^{\frac{1-\delta}{4}}g_N(z)\tilde{N}(dz,d\rho)\Big\|^2\Big]\Big\|_{L^p(\Omega,\R)} \\
  &\quad+C\Big\|\sup_{t\in[0,T]}\E^t\Big[\sup_{s\in[t,T]}\E\Big[\sup_{r\in[s,s']}\Big\|\int^r_s\int_\chi E(r-\rho)(-A)^{\frac{1-\delta}{4}}g_N(z)\tilde{N}(dz,d\rho)\Big\|^2\Big]\Big]\Big\|_{L^p(\Omega,\R)} \\
  &\quad+C\Big\|\sup_{t\in[0,T]}\E^t\Big[\sup_{s\in[t,T]}\Big\|\Big(\int_0^t+\int_t^s\Big)\int_\chi E(s-\rho)(-A)^{\frac{1-\delta}{4}}g_N(z)\tilde{N}(dz,d\rho)\Big\|^2\Big]\Big\|_{L^p(\Omega,\R)}\\
  &\quad+C\Big\|\sup_{t\in[0,T]}\E^t\Big[\sup_{s\in[t,T]}\E\Big[\Big\|\int_s^{s'}\int_\chi E(s'-\rho)(-A)^{\frac{1-\delta}{4}}g_N(z)\tilde{N}(dz,d\rho)\Big\|^2\Big]\Big]\Big\|_{L^p(\Omega,\R)}\\
  &\leq C_{p,T}+C\Big\|\sup_{t\in[0,T]}\Big\|\int_0^t\int_\chi E(t-\rho)(-A)^{\frac{1-\delta}{4}}g_N(z)\tilde{N}(dz,d\rho)\Big\|^2\Big\|_{L^p(\Omega,\R)}\\
  &\quad+C\Big\|\sup_{t\in[0,T]}\E\Big[\sup_{r\in[t,T]}\Big\|\int^r_t\int_\chi E(r-\rho)(-A)^{\frac{1-\delta}{4}}g_N(z)\tilde{N}(dz,d\rho)\Big\|^2\Big]\Big\|_{L^p(\Omega,\R)}\leq C_{p,T}.
\end{align*}
(iv) \textit{The proof of \eqref{ineq:WNC}.} By using
\eqref{ineq:WN-sup} in Lemma~\ref{lem:WN-sup}, the proof of \eqref{ineq:WNC} is similar to that of \eqref{ineq:NNC}, which is omitted. 

Combining (i)--(iv), the proof is completed.
\end{proof}

\begin{proof}[\textbf{Proof of Proposition \ref{prop:YN}}]
(\romannumeral1)
Owing to the linear growth condition of $F$, we derive that for any $t\in[0,T]$ and $p\ge 2$,
 \begin{align*}
&\quad \big\|(-A)^{-\frac{1-\delta}{4}}(Y^N(t)-Y^N(\lfloor t\rfloor))\big\|_{L^{p}(\Omega,H)}\\
&\leq \int_{\lfloor t\rfloor}^{t}\big\|F_N(X^N(s))\big\|_{L^{p}(\Omega,H)}ds+\big\|(-A)^{-\frac{1-\delta}{2}}(E(t-\lfloor t\rfloor)-\mathrm{I})\big\|_{\IL(H)}\times \\
  &\quad \int_0^{\lfloor t\rfloor}\big\|(-A)^{\frac{1-\delta}{4}}E(\lfloor t\rfloor-s)\big\|_{\IL(H)}\big\|F_N(X^N(s))\big\|_{L^{p}(\Omega,H)}ds
  \leq C_{p,T}(\Delta t)^{\frac{1-\delta}{2}}.
 \end{align*}	
 
(\romannumeral2)
According to the estimates obtained in Propositions~\ref{prop:NN}--\ref{prop:WNN-conditional}, we can derive that terms 
\begin{align*}
	&\big\|\sup_{t\in[0,T]}\E^t\big[\sup_{\zeta\in[0,T]}\big\|\IW^{N}(\zeta)\big\|^2\big]\big\|_{L^p(\Omega,\R)}, \qquad
	\big\|\sup_{t\in[0,T]}\E^t\big[\sup_{\zeta\in[0,T]}\big\|\mathcal N^N(\zeta)\big\|^2\big]\big\|_{L^p(\Omega,\R)},\\ 
	&\qquad\qquad\qquad\big\|\sup_{t\in[0,T]}\E^t\big[\sup_{s\in[t,T]}\E^{s}\big[\sup_{r\in[0,T]}\big\|\IW^{N}(r)\big\|^2\big] \big]\big\|_{L^p(\Omega,\R)},\\ 
&\qquad\qquad\qquad \big\|\sup_{t\in[0,T]}\E^t\big[\sup_{s\in[t,T]}\E^{s}\big[\sup_{r\in[0,T]}\big\|\mathcal N^N(r)\big\|^2\big] \big]\big\|_{L^p(\Omega,\R)},\\ 
  &\qquad\qquad\qquad\big\|\sup_{t\in[0,T]}\E^t\big[\sup_{s\in[t,T]}\E^{s}\big[\sup_{r\in[0,T]}\E^{s'}\big[\big\|\IW^{N}(r)\big\|^2\big]\big] \big]\big\|_{L^p(\Omega,\R)},\\ 
&\qquad\qquad\qquad \big\|\sup_{t\in[0,T]}\E^t\big[\sup_{s\in[t,T]}\E^{s}\big[\sup_{r\in[0,T]}\E^{s'}\big[\big\|\mathcal N^N(r)\big\|^2\big] \big]\big]\big\|_{L^p(\Omega,\R)}
\end{align*}
are all bounded for $p\ge 2$ with some constant $ C_{p,T}$. 

Recall that $Y^N$ satisfies the perturbed stochastic equation \eqref{eq:perturbed}. 
Using the linear growth condition of $F$, it implies that for all $\zeta\in[0,T]$,
\begin{align*}
&\quad\big\|\sup_{t\in[0,T]}\E^t\big[\|Y^N(\zeta)\|^2\big]\big\|_{L^p(\Omega,\R)}\\
&\leq  C\|x_0\|_{L^{2p}(\Omega,H)}^2\!+C_T\Big\|\sup_{t\in[0,T]}\E^t\Big[\int_0^{\zeta}\big\|F_N(Y^N(\eta)\!+\!\IW^{N}(\eta)\!+\!\mathcal N^N(\eta))\big\|^2d\eta\Big]\Big\|_{L^p(\Omega,\R)}\\
    &\leq C+C_T\Big\|\sup_{t\in[0,T]}\E^t\Big[\int_0^{\zeta}\|Y^N(\eta)\|^2d\eta\Big]\Big\|_{L^p(\Omega,\R)}\\
 &\quad+C_T\big\|\sup_{t\in[0,T]}\E^t\big[\sup_{\zeta\in[0,T]}\|\IW^{N}(\zeta)\|^2\big]\big\|_{L^p(\Omega,\R)}+C_T\big\|\sup_{t\in[0,T]}\E^t\big[\sup_{\zeta\in[0,T]}\|\mathcal N^N(\zeta)\|^2\big]\big\|_{L^p(\Omega,\R)}\\
 &\leq C_{p,T}+C_T\int_0^{\zeta}\big\|\sup_{t\in[0,T]}\E^t\big[\|Y^N(\eta)\|^2\big]\big\|_{L^p(\Omega,\R)}d\eta.
\end{align*}
Applying the Gr\"onwall inequality leads to 
$\sup\limits_{\zeta\in[0,T]}\big\|\sup\limits_{t\in[0,T]}\E^t\big[\big\|Y^N(\zeta)\big\|^2\big]\big\|_{L^p(\Omega,\R)}\leq C_{p,T}.$  

In addition,  for all $\varrho\in[0,T]$,
\begin{align*}
  &\quad\big\|\sup_{t\in[0,T]}\E^t\big[\sup_{s\in[t,T]}\E^{s}\big[\|Y^N(\varrho)\|^2\big]\big]\big\|_{L^p(\Omega,\R)}\\
  &\leq  C\|x_0\|_{L^{2p}(\Omega,H)}^2+C_T\Big\|\sup_{t\in[0,T]}\E^t\Big[\sup_{s\in[t,T]}\E^{s}\Big[\int_0^{\varrho}\big\|F_N(Y^N(r)+\IW^{N}(r)+\mathcal N^N(r))\big\|^2dr\Big] \Big]\Big\|_{L^p(\Omega,\R)}\\
    &\leq C_{T}+C_T\int_0^{\varrho}\big\|\sup_{t\in[0,T]}\E^t\big[\sup_{s\in[t,T]}\E^{s}\big[\|Y^N(r)\|^2\big]\big]\big\|_{L^p(\Omega,\R)}dr\\
    &\quad+C_T\big\|\sup_{t\in[0,T]}\E^t\big[\sup_{s\in[t,T]}\E^{s}\big[\sup_{r\in[0,T]}\|\IW^{N}(r)\|^2\big] \big]\big\|_{L^p(\Omega,\R)}\\
    &\quad+C_T\big\|\sup_{t\in[0,T]}\E^t\big[\sup_{s\in[t,T]}\E^{s}\big[\sup_{r\in[0,T]}\|\mathcal N^N(r)\|^2\big] \big]\big\|_{L^p(\Omega,\R)}\\
    &\leq C_{p,T}+C_T\int_0^{\varrho}\big\|\sup_{t\in[0,T]}\E^t\big[\sup_{s\in[t,T]}\E^{s}\big[\|Y^N(r)\|^2\big]\big]\big\|_{L^p(\Omega,\R)}dr.
\end{align*}
The Gr\"onwall inequality leads to 
$\sup\limits_{\varrho\in[0,T]}\big\|\sup\limits_{t\in[0,T]}\E^t\big[\sup\limits_{s\in[t,T]}\E^{s}\big[\|Y^N(\varrho)\|^2\big]\big]\big\|_{L^p(\Omega,\R)}\leq C_{p,T}.$ 

Combining the above estimates, we obtain 
\begin{align}\label{ineq:YN-estimate1}
    &\quad\Big\|\sup_{t\in[0,T]}\E^t\Big[\sup_{s\in[t,T]}\E^s\big[\big\|(-A)^{\frac{1-\delta}{4}}Y^N(\cdot)\big\|^2_{\mathcal{C}_2^0|\mathcal{F}_{s'},[s',T]}\big]\Big]\Big\|_{L^p(\Omega,\R)}\nonumber\\
    &\leq C\big\|(-A)^{\frac{1-\delta}{4}}x_0\big\|_{L^{2p}(\Omega,H)}^2\nonumber\\
    &\quad+\!C_T\Big\|\!\sup_{t\in[0,T]}\!\E^t\Big[\!\sup_{s\in[t,T]}\!\E^s\Big[\!\sup_{r\in[s',T]}\!\E^{s'}\Big[\int_0^r\!\big\|(-A)^{\frac{1\!-\!\delta}{4}}E(r\!-\!\rho)F(Y^N(\rho)\!+\!\IW^{N}(\rho)\!+\!\mathcal N^N(\rho))\big\|^2d\rho\Big]\Big]\Big]\Big\|_{L^p(\Omega,\R)}\nonumber\\
    &\leq C+C_T\Big\|\sup_{t\in[0,T]}\E^t\Big[\sup_{s\in[t,T]}\E^s\Big[\sup_{r\in[s',T]}\sup_{\rho\in[0,r]}\E^{s'}\big[\|Y^N(\rho)\|^2\big]\int_0^r\|E(r-\eta)(-A)^{\frac{1-\delta}{4}}\|^2_{\IL(H)}d\eta\Big]\Big]\Big\|_{L^p(\Omega,\R)}\nonumber\\
    &\quad+C_T\Big\|\sup_{t\in[0,T]}\E^t\big[\sup_{s\in[t,T]}\E^{s}\big[\sup_{r\in[0,T]}\E^{s'}\big[\big\|\IW^{N}(r)\big\|^2\big]\big] \big]\Big\|_{L^p(\Omega,\R)}\nonumber\\
    &\quad+C_T\Big\|\sup_{t\in[0,T]}\E^t\big[\sup_{s\in[t,T]}\E^{s}\big[\sup_{r\in[0,T]}\E^{s'}\big[\big\|\mathcal N^N(r)\big\|^2\big] \big]\big]\Big\|_{L^p(\Omega,\R)}\nonumber\\
    &\leq C_{p,T}+C_T\big\|\sup_{t\in[0,T]}\E^t\big[\sup_{s\in[t,T]}\E^s\big[\sup_{\rho\in[0,T]}\E^{s'}\big[\|Y^N(\rho)\|^2\big]\big]\big]\big\|_{L^p(\Omega,\R)}.
    \end{align}
 Note that 
    \begin{align}\label{ineq:YN-estimate2}
  &\quad \big\|\sup_{t\in[0,T]}\E^t\big[\sup_{s\in[t,T]}\E^s\big[\sup_{\rho\in[0,T]}\E^{s'}\big[\|Y^N(\rho)\|^2\big]\big]\big]\big\|_{L^p(\Omega,\R)}\nonumber\\
    &\leq  C_T\|x_0\|_{L^{2p}(\Omega,H)}^2\nonumber\\
    &\quad+C_T\Big\|\sup_{t\in[0,T]}\E^t\Big[\sup_{s\in[t,T]}\E^s\Big[\int_0^{T}\E^{s'}\big[\|Y^N(\vartheta)+\IW^N(\vartheta)+\mathcal N^N(\vartheta)\|^2\big]d\vartheta\Big]\Big]\Big\|_{L^p(\Omega,\R)}\nonumber\\
    &\leq C_{p,T}+C_T\int_0^{T}\big\|\sup_{t\in[0,T]}\E^t\big[\sup_{s\in[t,T]}\E^{s}\big[\|Y^N(\vartheta)+\IW^N(\vartheta)+\mathcal N^N(\vartheta)\|^2\big]\big]\big\|_{L^p(\Omega,\R)}d\vartheta\nonumber\\
    &\leq C_{p,T}+C_T\int_0^{T}\big\|\sup_{t\in[0,T]}\E^t\big[\sup_{s\in[t,T]}\E^{s}\big[\|Y^N(\vartheta)\|^2\big]\big]\big\|_{L^p(\Omega,\R)}d\vartheta.
\end{align}
Thus inserting \eqref{ineq:YN-estimate2} into \eqref{ineq:YN-estimate1} and then applying the Gr\"onwall inequality yield the desired result. The proof is completed.
\end{proof}

\section*{Appendix.}\label{sec6}

 \begin{proof}[\textbf{Proof of Proposition~\ref{prop:wellposedness}}]
 For brevity, we set $Y_1(t):=X(t)-\mathcal{W}(t)$ for all $t\in[0,T]$. Then $\{Y_1(t)\}_{t\in[0,T]}$ satisfies
\begin{equation*}
    \left\{
    \begin{aligned}
     &dY_1(t)=AY_1(t)dt+F\big(Y_1(t)+\mathcal{W}(t)\big)dt+\int_\chi G\big(Y_1(t)+\mathcal{W}(t), z\big)\tilde{N}(dz,dt), \quad t\in(0,T],  \\
     &Y_1(0)=x_0.
    \end{aligned}
    \right.
\end{equation*} 
For all $p\geq2$ and all $\alpha\in[0,\frac12)$, we note that there exists a constant $C_{p,T}>0$ such that $\sup_{t\in[0,T]}\|\IW(t)\|_{L^p(\Omega,\dot{H}^{\alpha})}\leq C_{p,T}.$ Using Assumption~\ref{assum:nonlinearity} and Lemma~\ref{ineq:maximalinequality} leads to 
 \begin{align*}
&\quad\|(-A)^{\frac{\alpha}{2}}Y_1(t)\|_{L^p(\Omega, H)}\\
&\leq \|(-A)^{\frac{\alpha}{2}}E(t)x_0\|_{L^p(\Omega, H)}+\int_0^t\big\|(-A)^{\frac{\alpha}{2}}E(t\!-\!s)F\big(Y_1(s)+\mathcal{W}(s)\big)\big\|_{L^p(\Omega, H)}ds\\
   &\quad+\Big\|\int_0^t\int_\chi (-A)^{\frac{\alpha}{2}}E(t\!-\!s)G\big(Y_1(s)+\mathcal{W}(s), z\big)\tilde{N}(dz,ds)\Big\|_{L^p(\Omega, H)}\\
   &\leq \|x_0\|_{L^p(\Omega, \dot{H}^\alpha)}+C\int_0^t\|(-A)^{\frac{\alpha}{2}}E(t\!-\!s)\|_{\IL(H)}\big(1+\|Y_1(s)\|_{L^p(\Omega, H)}+\|\mathcal{W}(s)\|_{L^p(\Omega, H)}\big)ds\\
   &\quad+C_{p,T}\Big(\int_0^t\E\Big[\int_\chi\big\|(-A)^{\frac{\alpha}{2}}G\big(Y_1(s)+\mathcal{W}(s),z\big)\big\|^p\nu(dz)\\
   &\quad+\Big(\int_\chi\big\|(-A)^{\frac{\alpha}{2}}G\big(Y_1(s)+\mathcal{W}(s), z\big)\big\|^2\nu(dz)\Big)^{\frac{p}{2}}\Big]ds\Big)^{\frac1p}.
   \end{align*}
   When $\alpha=0$, owing to the moment estimates of $\mathcal W$, Assumption~\ref{assum:jumpcoefficient}, and the H\"older inequality, we arrive at
\begin{align*}
   \|Y_1(t)\|^p_{L^p(\Omega, H)}&\leq \|x_0\|^p_{L^p(\Omega,H)}+C_{p,T}\int_0^t\big(1+\|Y_1(s)\|^p_{L^p(\Omega, H)}\big)ds\\
   &\quad+C_{p,T}\int_0^t\E\Big[\int_\chi\big(|g_1(z)|\big\|Y_1(s)+\mathcal{W}(s)\big\|+\|g(z)\|\big)^p\nu(dz)\\
   &\quad+\Big(\int_\chi\big(|g_1(z)|\big\|Y_1(s)+\mathcal{W}(s)\big\|+\|g(z)\|\big)^2\nu(dz)\Big)^{\frac{p}{2}}\Big]ds\\
   &\leq C_{p,T}+\|x_0\|^p_{L^p(\Omega,H)}+C_{p,T}\int_0^t\|Y_1(s)\|^{p}_{L^p(\Omega, H)}ds.
\end{align*}
Applying the Gr\"onwall inequality yields that $\|Y_1(t)\|_{L^p(\Omega, H)}\leq C_{p,T}$. 

Then for $\alpha\in(0,\frac12)$, we have
\begin{align*}
 \|(-A)^{\frac{\alpha}{2}}Y_1(t)\|^p_{L^p(\Omega, H)}&\leq C_p\|x_0\|^p_{L^p(\Omega,\dot{H}^\alpha)}+C_{p,T}\Big(\int_0^t(t\!-\!s)^{-\frac{\alpha}{2}}ds\Big)^p\sup_{s\in[0,t]}\big(1+\|Y_1(s)\|^p_{L^p(\Omega, H)}\big)\\
   &\quad+C_{p,T}\int_0^t\E\Big[\int_\chi\big\|g_{1}(z)(-A)^{\frac{\alpha}{2}}\big(Y_1(s)+\mathcal{W}(s)\big)+(-A)^{\frac{\alpha}{2}}g(z)\big\|^p\nu(dz)\\
   &\quad+\Big(\int_\chi\big\|g_{1}(z)(-A)^{\frac{\alpha}{2}}\big(Y_1(s)+\mathcal{W}(s)\big)+(-A)^{\frac{\alpha}{2}}g(z)\big\|^2\nu(dz)\Big)^{\frac{p}{2}}\Big]ds\\
   &\leq C_{p,T}+C_{p,T}\int_0^t\|(-A)^{\frac{\alpha}{2}}Y_1(s)\|^p_{L^p(\Omega, H)}ds.
 \end{align*}
Applying the Gr\"onwall inequality yields the desired assertion.
\end{proof}

\begin{proof}[\textbf{Proof of Lemma \ref{lem:WN-sup}}]
(\romannumeral1) 
   Denote $\IW^N_{s,t}:=\int_s^tE(t-r)P_NdW(r)$ for all $s\in[0,T]$ and $t\in[s,T]$. Based on the factorization method, we have expression $\mathcal W^N_{s,t}=\int_s^t(t-r)^{\alpha-1}e^{A(t-r)}\mathcal Z^N_{s,r}dr$ with $\mathcal Z^N_{s,r}=\frac{\sin (\pi\alpha)}{\alpha}\int_s^re^{A(r-\sigma)}(r-\sigma)^{-\alpha}P_N dW(\sigma)$ for $r\in[s,t]$; see e.g. \cite[Proposition 5.9, Theorem 5.10]{GDP-2014} for details. Note that for all $k\geq1$,
\begin{align}\label{fact1}
& \mathbb E\Big[\sup_{t\in[s,T]}\|(-A)^{\frac{1-\delta}{4}}\mathcal W^N_{s,t}\|^{2k}\Big]
\leq \mathbb E\Big[\sup_{t\in[s,T]}\Big|\int_s^t(t-\rho)^{\alpha-1}e^{-\lambda_1(t-\rho)}\|(-A)^{\frac{1-\delta}{4}}\mathcal Z^N_{s,\rho}\| d\rho\Big|^{2k}\Big]\notag\\
&\leq C\mathbb E\Big[\sup_{t\in[s,T]}\Big(\int_s^t(t-\rho)^{p(\alpha-1)}e^{-p\lambda_1(t-\rho)}d\rho\Big)^{\frac{2k}{p}}\Big(\int_s^T\|(-A)^{\frac{1-\delta}{4}}\mathcal Z^N_{s,\rho}\|^q d\rho\Big)^{\frac{2k}{q}}\Big],\nonumber
\end{align}
where we utilized the H\"older inequality with $p,q\geq1$, $\frac1p+\frac1q=1$ in the last step. When $2k\ge q$, by the H\"older inequality and the Burkholder--Davis--Gundy inequality, we derive
\begin{align*}
&\quad \mathbb E\Big[\Big(\int_s^T\|(-A)^{\frac{1-\delta}{4}}\mathcal Z^N_{s,r}\|^q dr\Big)^{\frac{2k}{q}}\Big]\leq C_T\mathbb E\Big[\int_s^T\|(-A)^{\frac{1-\delta}{4}}\mathcal Z^N_{s,r}\|^{2k} dr\Big]\\
&\leq C_T\int_s^T\Big(\int_s^r(\frac{\sin(\pi \alpha)}{\alpha})^2(r-\sigma)^{-2\alpha}\|(-A)^{\frac{2-\delta+\epsilon}{4}}e^{A(r-\sigma)}\|^2_{\mathcal L(H)}\|(-A)^{-\frac{1+\epsilon}{4}}\|^2_{\mathcal L_2(H)} d\sigma\Big)^k dr\\
&\leq C_T\int_s^T\Big(\int_s^r(r-\sigma)^{-2\alpha-1+\frac{\delta-\epsilon}{2}}d\sigma\Big)^k dr,
\end{align*}
which is finite if $-2\alpha+\frac{\delta-\epsilon}{2}>0,$ i.e., $\alpha<\frac{\delta-\epsilon}{4}.$ It means that $\mathcal Z^N_{s,\cdot}\in L^{2k}(\Omega\times[s,T],D((-A)^{\frac{1-\delta}{4}}))$. 
Hence, when positive parameters $\alpha$ and $p$ are chosen such that  $\alpha<\frac{\delta-\epsilon}{4}$ and $p<\frac{1}{1-\alpha}$, we derive that 
$
\sup_{t\in[s,T]}\Big(\int_s^t(t-\rho)^{p(\alpha-1)}e^{-p\lambda_1(t-\rho)}\mathrm d\rho\Big)^{\frac{2k}{p}}\leq C_{k,p,T}. 
$
Therefore, for all $k\geq k_0$ with large number $k_0\geq\frac{q}{2}$, we obtain 
$
    \mathbb E\Big[\sup_{t\in[s,T]}\|(-A)^{\frac{1-\delta}{4}}\mathcal W^N_{s,t}\|^{2k}\Big]<\infty.
$
For $k<k_0$, we have
\begin{align*}
    \mathbb E\Big[\sup_{t\in[s,T]}\|(-A)^{\frac{1-\delta}{4}}\mathcal W^N_{s,t}\|^{2k}\Big]\leq \mathbb E\Big[\sup_{t\in[s,T]}\|(-A)^{\frac{1-\delta}{4}}\mathcal W^N_{s,t}\|^{2k_0}\Big]<\infty. 
\end{align*}

(\romannumeral2) 
For $p\geq1$ and any $t\in[0,T]$, applying the Burkholder--Davis--Gundy inequality and properties~\eqref{prop:semigroup1}--\eqref{prop:semigroup2} gives 
\begin{align*}
&~~~\big\|(-A)^{-\frac{1-\delta}{4}}\big(\IW^N(t)-\IW^N(\lfloor t\rfloor)\big)\big\|_{L^{2p}(\Omega,H)}\leq \Big\|\int_{\lfloor t\rfloor}^t (-A)^{-\frac{1-\delta}{4}}E(t-r)P_NdW(r)\Big\|_{L^{2p}(\Omega,H)}\\
&\quad+\Big\|\int_0^{\lfloor t\rfloor} (-A)^{-\frac{1-\delta}{4}}(E(t\!-\!\lfloor t\rfloor)\!-\!\mathrm{I})E(\lfloor t\rfloor\!-\!r)P_NdW(r)\Big\|_{L^{2p}(\Omega,H)}\\
 &\leq C_p\Big[\Big|\int_{\lfloor t\rfloor}^t(t-r)^{-\delta}\|(-A)^{-\frac{1+\delta}{4}}\|_{\IL_2(H)}^2dr\Big|^{\frac{1}{2}}+(\Delta t)^{\frac{1-\delta}{2}}\Big|\int_{0}^{\lfloor t\rfloor}(\lfloor t\rfloor-r)^{-\frac{1-\delta}{2}}e^{-\lambda_1(\lfloor t\rfloor-r)}dr\Big|^{\frac{1}{2}}\Big]\\
 &\leq C_{p,T}(\Delta t)^{\frac{1-\delta}{2}},
\end{align*}
where we used $\int_0^\infty x^{\eta-1}e^{-x}dx<\infty, \,\eta>0$. 
The proof is completed.
\end{proof}

\bibliographystyle{plain}
\bibliography{references.bib}

 \end{document}